\documentclass[11pt, a4paper]{amsart}
\usepackage[T1]{fontenc}
\usepackage{stix2}

\usepackage{geometry}
\usepackage{microtype} 

\DeclareRobustCommand{\SkipTocEntry}[5]{}

\usepackage{amsthm}
\usepackage{mathtools} 

\usepackage{tikz-cd}
\tikzcdset{arrow style=tikz, diagrams={>=stealth}} 

\usepackage{graphicx}
\usepackage{rotating} 
\usepackage[pdfpagelabels]{hyperref}
\hypersetup{
  pdftitle={The profinite shape equivalence},
  pdfauthor={Tim Holzschuh},
  pdfsubject={},
  pdfkeywords={},
  colorlinks=true,    
  linkcolor=blue,     
  citecolor=blue,     
  filecolor=blue,      
  urlcolor=blue,       
  breaklinks=true,
  bookmarksopen=true,
  bookmarksnumbered=true,
  pdfpagemode=UseOutlines,
  plainpages=false,
  unicode=true
  }
  
\numberwithin{equation}{subsection}

\usepackage{color}
\usepackage{soul}
\setuldepth{Berlin}

\usepackage{csquotes}

\usepackage[nameinlink]{cleveref} 

\creflabelformat{theorem}{#2\textup{#1}#3} 
\creflabelformat{lemma}{#2\textup{#1}#3}
\creflabelformat{equation}{#2(#1)#3} 
\crefformat{section}{\S#2#1#3}
\crefformat{subsection}{\S#2#1#3}
\crefformat{subsubsection}{\S#2#1#3}

\crefname{recollection}{Recollection}{Recollections}
\crefname{definition}{Definition}{Definitions}
\crefname{thmABC}{Theorem}{Theorems}

\newtheorem{theorem}{Theorem}[section]
\newtheorem{lemma}[theorem]{Lemma}
\newtheorem{proposition}[theorem]{Proposition}
\newtheorem{corollary}[theorem]{Corollary}

\newtheorem{thmABC}{Theorem}

\theoremstyle{definition}
\newtheorem{definition}[theorem]{Definition}

\newtheorem{remark}[theorem]{Remark}
\newtheorem{recollection}[theorem]{Recollection}
\newtheorem{notation}[theorem]{Notation}

\newtheorem*{definition*}{Definition}
\newtheorem*{remark*}{Remark}

\usepackage[shortlabels]{enumitem}
\setlist{leftmargin=.8cm}

\makeatletter
\def\cref@thmoptarg[#1]#2#3#4{%
    \ifhmode\unskip\unskip\par\fi%
    \normalfont%
    \trivlist%
    \let\thmheadnl\relax%
    \let\thm@swap\@gobble%
    \thm@notefont{\fontseries\mddefault\upshape}%
    \thm@headpunct{.}
    \thm@headsep 5\p@ plus\p@ minus\p@\relax%
    \thm@space@setup%
    #2
    \@topsep \thm@preskip               
    \@topsepadd \thm@postskip           
    \def\@tempa{#3}\ifx\@empty\@tempa%
      \def\@tempa{\@oparg{\@begintheorem{#4}{}}[]}%
    \else%
      \refstepcounter[#1]{#3}
      \@namedef{cref@#3@alias}{#1}
      \def\@tempa{\@oparg{\@begintheorem{#4}{\csname the#3\endcsname}}[]}%
    \fi%
    \@tempa}%
\makeatother

\usepackage[style=alphabetic, backend=biber, natbib=true,
    url=false,
    doi=false, 
    isbn=false,
    eprint=false, maxbibnames=99]{biblatex}
\DeclareFieldFormat[article, inbook, incollection, inproceedings, misc, thesis, unpublished]{title}{\textit{#1}}
\DeclareFieldFormat{journaltitle}{{\rmfamily #1}}
\DeclareFieldFormat{booktitle}{{\rmfamily #1}}
\DeclareFieldFormat{pages}{#1}
\renewrobustcmd*{\bibinitdelim}{\addnbthinspace{}}
\renewrobustcmd*{\bibnamedelima}{\addnbthinspace{}}
\renewrobustcmd*{\bibnamedelimd}{\addnbthinspace{}}

\DeclareBibliographyDriver{article}{%
  \usebibmacro{author/editor}%
  \newunit
  \usebibmacro{title}%
  \newunit
  \usebibmacro{journal}%
  \setunit*{\addspace}%
  \printtext{\textbf{\printfield{volume}}}%
  \setunit*{\addspace}%
  \printtext[parens]{\printfield{year}}%
  \setunit{\addcomma\space}%
  \iffieldundef{number}{}{\printtext{no.~\printfield{number}}}%
  \setunit{\addcomma\addspace}%
  \printfield{pages}%
  \iffieldundef{pubstate}{}{\addcomma\addspace\printfield{pubstate}%
  	 \iffieldundef{doi}{}{\addcomma\addspace\printfield{doi}}%
	 }%
  \usebibmacro{finentry}}

\DeclareSourcemap{
  \maps[datatype=bibtex]{
    \map{
      \step[fieldset=urldate,  null]
      \step[fieldset=pubstate, null]
      \step[fieldset=note,     null]   
    }
  }
}

\newcommand{\stackstag}[1]{\href{https://stacks.math.columbia.edu/tag/#1}{Tag #1}}

\newcommand{\kerodontag}[1]{\href{https://kerodon.net/tag/#1}{Tag #1}}

\newcommand{\SAGapp}[1]{\href{https://www.math.ias.edu/~lurie/papers/SAG-rootfile.pdf\#appendix.#1}{Appendix #1}}

\newcommand{\HTTthm}[2]{\href{https://www.math.ias.edu/~lurie/papers/HTT.pdf\#theorem.#2}{#1 #2}}

\newcommand{\SAGthm}[2]{\href{https://www.math.ias.edu/~lurie/papers/SAG-rootfile.pdf\#theorem.#2}{#1 #2}}

\newcommand{\stacks}[1]{\cite[\stackstag{#1}]{stacks-project}}
\newcommand{\kerodon}[1]{\cite[\kerodontag{#1}]{kerodon}}
\newcommand{\HTT}[2]{\cite[\HTTthm{#1}{#2}]{HTT}}

\newcommand{\SAG}[2]{\cite[\SAGthm{#1}{#2}]{SAG}}

 \newcommand{\from}{\colon}
 \newcommand{\defeqq}{\coloneq}

 \newcommand{\lref}[1]{\ref{#1}}

 \newcommand{\ot}{\leftarrow}
 \newcommand{\epi}{\twoheadrightarrow}
 \newcommand{\surj}{\epi}
 \newcommand{\mono}{\hookrightarrow}

 \DeclarePairedDelimiter{\sqrbr}{[}{]}

 \DeclarePairedDelimiterX{\setcond}[2]{\{}{\}}{#1 \;\delimsize\vert\; #2}

 \newcommand{\blank}{-}

 \DeclareMathOperator{\im}{im}

\newcommand{\isomto}{\xlongrightarrow{\,\smash{\raisebox{-0.95ex}{\ensuremath{\displaystyle\sim}}}\,}}

\DeclareMathOperator{\id}{id}

 \newcommand{\ladj}{\dashv}

\renewcommand{\in}{\smallin}

\renewcommand{\ast}{{\mathbin{*}}}

 \newcommand{\et}{\operatorname{\acute{e}t}}

 \newcommand{\cts}{\operatorname{cts}}

 \newcommand{\fin}{\operatorname{fin}}

 \newcommand{\bN}{{\mathbb N}}

 \newcommand{\cC}{{\mathscr C}}
 
 \newcommand{\cE}{{\mathscr E}}

 \newcommand{\cS}{{\mathscr S}}

 \newcommand{\cX}{{\mathscr X}}
 \newcommand{\cY}{{\mathscr Y}}

 \newcommand{\cat}[1]{\mathscr{#1}}
 \newcommand{\catC}{\cat{C}}
 \newcommand{\catD}{\cat{D}}
 \newcommand{\catE}{\cat{E}}

 \newcommand{\op}{\operatorname{op}}

 \DeclareMathOperator{\Set}{\mathbf{Set}}
 \DeclareMathOperator{\Grp}{\mathbf{Grp}}
 \DeclareMathOperator{\Grpd}{\mathbf{Grpd}}
 
 \DeclareMathOperator{\Ab}{\mathbf{Ab}}

 \DeclareMathOperator{\Ani}{\mathbf{Ani}}
 \DeclareMathOperator{\Cat}{\mathbf{Cat}}

 \newcommand{\map}{\Map}
 \DeclareMathOperator{\Map}{Map}

 \DeclareMathOperator{\Fun}{Fun}

 \DeclareMathOperator{\colim}{colim}
 \DeclareMathOperator{\limit}{lim}

 \DeclareMathOperator{\Pro}{Pro}
 \DeclareMathOperator{\Ind}{Ind}

 \DeclareSymbolFont{supplsymbols}{T1}{\familydefault}{m}{n}
 \SetSymbolFont{supplsymbols}{bold}{T1}{\familydefault}{bx}{n}

 \DeclareMathSymbol{\frenchlq}{\mathopen}{supplsymbols}{19}
 \DeclareMathSymbol{\frenchrq}{\mathclose}{supplsymbols}{20}

 \makeatletter
   \def\varphlim@#1#2{%
     \vtop{\m@th\ialign{##\cr
       \hfil$#1\operator@font\phantom{lim}$\hfil\cr
       \noalign{\nointerlineskip\kern1.5\ex@}#2\cr
       \noalign{\nointerlineskip\kern-\ex@}\cr}}%
   }
   \def\varinjcolim{%
     \mathop{\ooalign{colim\cr\hfil$\mathpalette\varphlim@{\rightarrowfill@\textstyle}$\hfil}}\nmlimits@%
   }
 \makeatother

 \DeclareMathOperator{\prolimit}{\frenchlq \limit \frenchrq}

 \newcommand{\point}{\ast}

 \newcommand{\rel}{\!\downarrow\!}

 \newcommand{\cartesian}{\arrow[dr, phantom, very near start, "{ \lrcorner }"]}

\usetikzlibrary{decorations.markings}

\makeatletter
\tikzcdset{
  open/.code     = {\tikzcdset{hook, circled};},
  closed/.code   = {\tikzcdset{hook, slashed};},
  open'/.code    = {\tikzcdset{hook', circled};},
  closed'/.code  = {\tikzcdset{hook', slashed};},
  circled/.code  = {\tikzcdset{markwith = {\draw (0,0) circle (.375ex);}};},
  slashed/.code  = {\tikzcdset{markwith = {\draw[-] (-.4ex,-.4ex) -- (.4ex,.4ex);}};},
  markwith/.code ={
    \pgfutil@ifundefined%
    {tikz@library@decorations.markings@loaded}%
    {\pgfutil@packageerror{tikz-cd}{You need to say %
      \string\usetikzlibrary{decorations.markings} to use arrows with markings}{}}{}%
    \pgfkeysalso{/tikz/postaction = {
      /tikz/decorate,
      /tikz/decoration={markings, mark = at position 0.5 with {#1}}}
    }
  },
}
\makeatother

 \DeclareMathOperator{\htpygrp}{\operatorname{\pi}}

 \DeclareMathOperator{\fib}{fib}

 \DeclareMathOperator{\lcc}{\operatorname{lcc}}

 \DeclareMathOperator{\Mon}{Mon}

 \DeclareMathOperator{\pfhtpytype}{\widehat{\Pi}_{\infty}}
 \NewDocumentCommand{\relpfhtpytype}{m}{\widehat{\Pi}_{\infty/ #1}}
 \DeclareMathOperator{\ethtpytype}{\widehat{\Pi}^{\et}_{\infty}}
 \NewDocumentCommand{\truncethtpytype}{m}{\widehat{\Pi}^{\et}_{#1}}
 \DeclareMathOperator{\etfdtlgrp}{\hat{\htpygrp}_{1}^{\et}}

 \newcommand{\pfcompl}{^{\wedge}_{\pi}}
 \NewDocumentCommand{\lcccompl}{o}{^{\wedge}_{\lcc\IfValueT{#1}{, #1}}}
 \newcommand{\ppfcompl}[1]{^{\wedge}_{#1}}
 \newcommand{\Spfcompl}{\ppfcompl{\Sigma}}

 \DeclareMathOperator{\ev}{ev}
 \DeclareMathOperator{\const}{const}

\DeclareMathOperator{\Str}{\partial}
\DeclareMathOperator{\Unstr}{\int}

\title{The profinite shape equivalence}

\author{Tim Holzschuh}
\address{Tim Holzschuh, Institut f\"{u}r Mathematik, Universit\"at Heidelberg \\ \newline \indent
Im Neuenheimer Feld 205, 69120 Heidelberg, Germany}
\email{tholzschuh@mathi.uni-heidelberg.de}

\date{\today}

\begin{document}

\begin{abstract}
  We show that, for every $\infty$-topos $\cX$\!, the profinite shape induces a canonical equivalence between the $\infty$-category of pro-objects in locally constant constructible sheaves of anima in $\cX$ and the slice $\infty$-category of profinite anima over the profinite shape $\pfhtpytype(\cX)$.
  We deduce this profinite shape equivalence from the monodromy classification theorem by establishing a continuous version of the straightening--unstraightening equivalence for profinite anima.
  In addition, we explain how the monodromy classification theorem recovers its classical counterparts concerning locally constant sheaves of finite sets and of (abelian) groups.
\end{abstract}

\maketitle

\setcounter{tocdepth}{1}  
{\scriptsize \tableofcontents}

\setcounter{section}{-1}

\section{Introduction}
\label{sec:introduction}
Profinite shape theory attaches to an $\infty$-topos $\cX$ a profinite anima $\pfhtpytype(\cX)$ distilling all of the cohomological information of locally constant constructible sheaves in $\cX$ into a homotopy type.
For example, for a scheme $S$, the profinite shape of its étale $\infty$-topos $S_{\et}$ is the profinite \emph{étale homotopy type} of $S$: a invariant in arithmetic geometry that has been used by, e.g., Schmidt and Stix in \cite{schmidt2016} to resolve one of the main conjectures of Grothendieck's anabelian programme.

Explicitly computing and controlling $\pfhtpytype(\cX)$ is generally very difficult.
One reason for this is that many techniques from classical homotopy theory are not available --- for instance, explicitly manipulating paths.
Instead, a central tool is the \emph{monodromy classification theorem}, which asserts a canonical equivalence
\[
    \Fun^{\cts}(\pfhtpytype(\cX), \Ani_{\pi}) \simeq \lcc(\cX)
\]
between the $\infty$-categories of continuous local systems of $\pi$-finite anima on $\pfhtpytype(\cX)$ and locally constant constructible sheaves in $\cX$ (see \Cref{thm:monodromy-equivalence} for a more precise formulation and \cite[(10.1)]{bachmann2020} as well as \cite[(4.4.18)]{barwick2020} for proofs).
This equivalence often allows to translate problems concerning $\pfhtpytype(\cX)$ into sheaf-theory.

The purpose of this note is to prove the following consequence of the monodromy classification theorem:

\begin{thmABC}[{profinite shape equivalence}]
  \label{thm:A}
  Let $\cX$ be an $\infty$-topos.
  The profinite shape functor induces an equivalence of $\infty$-categories
  \[
    \pfhtpytype \from \Pro(\lcc(\cX)) \to \Pro(\Ani_{\pi})_{/\pfhtpytype(\cX)}.
  \]
\end{thmABC}
See \Cref{thm:lcc-sheaves-and-profinite-shape} for a more precise
formulation that also spells out an explicit inverse as well as its functoriality under geometric morphisms, and \Cref{thm:pro-Sigma-shape-equivalence} for a variant for the $\Sigma$-completed shape functor for some full class of finite groups $\Sigma$.

\Cref{thm:A} is a formal consequence of the monodromy classification theorem combined with our second main result, a continuous analogue of the straightening--unstraightening equivalence for profinite anima:

\begin{thmABC}[continuous (un-)straightening]
  \label{thm:B}
  Let $B$ be a profinite anima.
  There is a canonical adjoint equivalence of $\infty$-categories
  \[
    \begin{tikzcd}
      \Pro(\Fun^{\cts}(B, \Ani_{\pi})) \arrow[r, shift left = 0.2em, "{\Unstr^{B}}"] & \Pro(\Ani_{\pi})_{/B}. \arrow[l, shift left = 0.2em, "{\Str_{B}}"]
    \end{tikzcd}
  \]
\end{thmABC}
Again, we give a more general pro-$\Sigma$ variant and also describe the functors $\Unstr^B \ladj \Str_{B}$ in \Cref{thm:continuous-un-straightening}.

To deduce \Cref{thm:A} from \Cref{thm:B}, one needs to unpack the monodromy equivalence explicitly.
Since in the existing literature, the functors furnishing the monodromy equivalence are left mostly implicit, we take this opportunity to explicate them.
To round things up, we also discuss how the monodromy classification theorem relates to more classical notions of local systems of finite sets and (abelian) groups considered in e.g. \cite{sga1}.

\Cref{thm:A} constitutes the main formal input in forthcoming joint work with Schmidt and Stix generalising Friedlander's fibre sequence in étale homotopy theory of \cite{friedlander1973a}.

\subsection*{Linear overview}
In \Cref{sec:background}, we recall some background on profinite anima and filtered colimits of $\infty$-categories; the familiar reader can safely skip this section.
Next we study continuous local systems on profinite anima in \Cref{sec:continuous-local-systems}.
We then prove \Cref{thm:B}, the continuous analogue of the straightening-unstraightening equivalence for profinite anima, in \Cref{sec:continous-un-straightening}.
Finally, in \Cref{sec:profinite-shape-theory}, we deduce the profinite shape equivalence, \Cref{thm:A}, and discuss some of its consequences.

\subsection*{Acknowledgments}
The author acknowledges support by Deutsche Forschungsgemeinschaft (DFG) through the Walter Benjamin Programme, project number 565137331, as well as the Collaborative Research Centre TRR 326 “Geometry and Arithmetic of Uniformized Structures”, project number 444845124.
We would also like to thank Alexander Schmidt and Jakob Stix for insisting that we should write up these notes.
Moreover, we thank Dustin Clausen for telling us about a useful result in \cite{clausen2025}, Peter Haine for various comments on an early draft, and IHÉS for providing excellent working conditions.


\section{Background}
\label{sec:background}
In this section, we will collect some background on profinite anima and filtered colimits of $\infty$-categories.
\subsection{Profinite anima}
\label{subsec:profinite-anima}

\begin{notation}[pro-categories]
  \label{not:pro-categories}
	Given an $\infty$-category $ \catC $, we write $ \Pro(\catC) $ for the $\infty$-category of pro-objects in $ \catC $ obtained by formally adjoining cofiltered limits to $ \catC $.
	The existence of $ \Pro(\catC) $ is a special case of (the dual of) \HTT{}{5.3.6.2}.
	Given a functor $ F \from \catC \to \catD $, we simply write $ F \from \Pro(\catC) \to \Pro(\catD)$ for the cofiltered-limit-preserving extension of $ F $.
\end{notation}

Recall that a \emph{full class of finite groups} is a collection of finite groups that is closed under taking subgroups, homomorphic images, and extensions.

\begin{definition}[(pro-)finite anima]
\label{def:profinite-anima}
Let $B$ be an anima and $\Sigma$ a nonempty full class of finite groups.
\begin{enumerate}[(1)]
	\item \label{defitem:pi-finite} 
        The anima $B$ is said to be \emph{$\Sigma$-finite} if it is truncated, has finitely many connected components, and if for all points $b \in B$ and all $i \geq 1$ the homotopy group $\htpygrp_i(B, b)$ is in $\Sigma$.
	
  	\item 
	We write $\Ani_{\Sigma} \subset \Ani$ for the full subcategory spanned by $\Sigma$-finite anima.
	
  	\item 
	We call $\Pro(\Ani_{\Sigma})$ the $\infty$-category of \emph{pro-$\Sigma$ anima}.

        \item
        We write $\pi$ for the class consisting of all finite groups and call $\Pro(\Ani_{\pi})$ the $\infty$-category of \emph{profinite anima}.
\end{enumerate}
\end{definition}

\begin{recollection}[profinite \& pro-$\Sigma$ completion]
  \label{rec:Spfcompl-preserves-cofiltered-limits}
  Let $\Sigma$ denote a nonempty full class of finite groups.
  \begin{enumerate}[(1)]
  \item \label{recitem:pfcompl}
  The canonical inclusions $\Pro(\Ani_{\Sigma}) \subset \Pro(\Ani_{\pi}) \subset \Pro(\Ani)$ admit left adjoints
  \[
    \Pro(\Ani) \xrightarrow{(\blank)\pfcompl}  \Pro(\Ani_{\pi}) \xrightarrow{(\blank)\Spfcompl}  \Pro(\Ani_{\Sigma})
  \]
  called \emph{profinite completion} and \emph{pro-$\Sigma$ completion}, respectively.
  \item \label{recitem:pfcompl-is-locally-cartesian-and-preserves-limits}
  The localisations 
  $(\blank)\pfcompl$ and $(\blank)\Spfcompl$ 
  preserve cofiltered limits and are locally cartesian, see \cite[Prop.~3.18]{haine2024b}.  
  The latter means that, for every map $K \to B$ of pro-$\Sigma$ anima and every map $E \to B$ of pro-anima, the natural map 
\[
\big(K \times_B E\big)\Spfcompl \isomto K \times_B E\Spfcompl
\]
is an equivalence (see \cites[\S 3]{haine2024b}[\S1.2]{gepner2017}[\S3.2]{hoyois2017a}).
  \end{enumerate}
\end{recollection}

An important feature of profinite anima is that equivalences can be detected on homotopy groups.

\begin{theorem}[profinite Whitehead theorem, {\SAG{}{E.3.1.6}}]
\label{thm:profinite-whitehead}
The following are equivalent for a map $f \from E \to B$ of profinite anima.
\begin{enumerate}[(1)]
	\item 
	The map $f$ is an equivalence of profinite anima.
	
    	\item 
	The induced map on connected components $f_{\! \ast} \from \htpygrp_{0}\!E \to \htpygrp_0\!B$ 
	is bijective and for every point $e \in E$ with image $b = f(e) \in B$ and every 
	$k \geq 1$ the induced homomorphism
    	\[
      		f_{\!\ast} \from \htpygrp_k(E, e) \to \htpygrp_k(B, b)
   	 \]
    	is an isomorphism of profinite groups.
\end{enumerate}
\end{theorem}

As a consequence, we have the following useful criterion for detecting pullback squares of profinite anima.

\begin{corollary}[fibrewise criterion of pullbacks]
    \label{cor:fibrewise-criterion-pullbacks}
    The following are equivalent for a commutative square
    \[
        \begin{tikzcd}
            E' \arrow[r] \arrow[d]
              & E \arrow[d] \\
            B' \arrow[r]
              & B
        \end{tikzcd}
    \]
    of  anima or of profinite anima, respectively:
    \begin{enumerate}[(1)]
        \item
        \label{lemitem:fibre product of anima}
        The square is cartesian.
        \item
        \label{lemitem:fibre of canonical map}
        For every point $p \in  E \times_{B} B'$, the fibre $\fib_p(E' \to E \times_{B} B')$ is contractible.
        \item
        \label{lemitem:homotopy equivalent fibres}
        For every point $b' \in B'$ with image $b \in B$, the induced map
            \[
                \fib_{b'}(E'\to B') \longrightarrow \fib_{b}(E\to B)
            \]
            is an equivalence.
    \end{enumerate}
\end{corollary}

\begin{proof}
  If \labelcref{lemitem:fibre of canonical map} holds, then the homotopy sequence for $E' \to E \times_{B} B'$ and every point $e' \in E'$ with image $p$ in $E \times_{B} B'$ yields isomorphisms $\pi_i(E',e') \to \pi_i(E \times_{B} B',p)$ and all $i \in \bN$, so that \labelcref{lemitem:fibre product of anima} follows from \Cref{thm:profinite-whitehead}.
  The implication \labelcref{lemitem:fibre product of anima} $\Rightarrow$ \labelcref{lemitem:homotopy equivalent fibres} is formal.

For the remaining direction \labelcref{lemitem:homotopy equivalent fibres} $\Rightarrow$
\labelcref{lemitem:fibre of canonical map}, we consider an arbitrary point $p \in E \times_{B} B'$ and its images $e\in E$, $b' \in B'$ and $b \in B$. A diagram chase in
\[
	\begin{tikzcd}
        \fib_{b'}(E'\to B')  
        \cartesian \arrow[d] \arrow[r]
        & E' \arrow[d]
        \\
	\fib_{b}(E\to B) \cartesian \arrow[d] \arrow[r]
        & E \times_{B} B' \cartesian \arrow[d] \arrow[r]
        & E \arrow[d]
        \\
        \point \arrow[r, "b'"]
        & B' \arrow[r] & B   .
     	\end{tikzcd}
\]
shows that $\fib_p(E'\to E\times_B B')$ is equivalent to the fibre of $\fib_{b'}(E'\to B') \to \fib_{b}(E\to B)$ in the point induced by $e$.
Since the fibre of an equivalence is contractible, this proves the claim.
\end{proof}

\begin{lemma}
  \label{lem:level-representation-profinite-anima}
  Let $B$ be a profinite anima.
  There exists a level-representation $B = \prolimit_{\lambda \in \Lambda} B_{\lambda}$ with the property that any of the projection maps $B \to B_{\mu}$ is surjective on connected components.
\end{lemma}

\begin{proof}
  Choose any level-representation $B = \prolimit_{\lambda \in \Lambda} B_{\lambda}$.
  Given any $\mu \in \Lambda$, write $p_{\mu} \from B \to B_{\mu}$ for the projection and $T_{\!\mu} \defeqq \im(\htpygrp_0(p_{\mu})) \subset \htpygrp_{0}(B_{\mu})$.
  The limit over $\mu \in \Lambda$ of the pullback diagrams
  \[
    \begin{tikzcd}
      B'_{\mu} \arrow[r] \arrow[d] \cartesian & B_{\mu} \arrow[d] \\
      T_{\!\mu} \arrow[r, hook] & \htpygrp_{0}(B_{\mu})
    \end{tikzcd}
  \]
 yields a pullback diagram
  \[
    \begin{tikzcd}
      \prolimit_{\lambda \in \Lambda} B'_{\lambda} \arrow[r] \arrow[d] \cartesian & B \arrow[d] \\
      \prolimit_{\lambda \in \Lambda} T_{\!\lambda} \arrow[r] & \htpygrp_0(B)
    \end{tikzcd}
  \]
  in $\Pro(\Ani_{\pi})$.
  By construction, the map $\prolimit_{\lambda \in \Lambda} T_{\!\lambda} \to \htpygrp_0(B)$
  is an equivalence.
  Hence so is $\prolimit_{\lambda \in \Lambda} B'_{\lambda} \to B$, and moreover, for each $\mu \in \Lambda$ the projection $\prolimit_{\lambda \in \Lambda} B'_{\lambda} \to B'_{\mu}$ is surjective on connected components.
\end{proof}

\subsection{Filtered colimits of \texorpdfstring{$\infty$}{infinity}-categories}
\label{subsec:filtered-colimits-of-infinity-categories}

\begin{recollection}[filtered colimits of $\infty$-categories]
  \label{rec:colimits-of-infinity-categories}
  Let $\Lambda$ be a filtered $\infty$-category. and $\catC_{\bullet} \from \Lambda \to \Cat_{\infty}$ a diagram.
  \begin{enumerate}[(1)]
    \item \label{recitem:colimit-cofinal} For any $\lambda \in \Lambda$, the forgetful functor $\Lambda_{\lambda/} \to \Lambda$ is colimit-cofinal, see \kerodon{06D0} and the dual of \kerodon{02XR}.
    \item\label{recitem:cocomplete} The $\infty$-category $\Cat_{\infty}$ of (small) $\infty$-categories admits all (small) limits and colimits. 
    \item \label{recitem:finite} An $\infty$-category $\catC$ is \emph{finite} if it is contained in the smallest full subcategory of $\Cat_{\infty}$ containing $\emptyset, \Delta^{0}$, and $\Delta^{1}$ that is closed under pushouts.
    \item \label{recitem:compact-categories} The compact objects in $\Cat_{\infty}$ are precisely the retracts of finite $\infty$-categories, see \cites[\kerodontag{068L}, \kerodontag{068E}]{kerodon}.
      In particular, finite $\infty$-categories are compact in $\Cat_{\infty}$.
    \item \label{recitem:finite-diagrams-lift} As a consequence of \lref{recitem:compact-categories}, any finite diagram $K \to \colim_{\lambda \in \Lambda} \catC_{\lambda}$ lifts to some finite stage $\mu \in \Lambda$.
    \item \label{recitem:maps-in-filtered-colimits} Given two objects $c = [c_{\lambda}]$ and $c' = [c'_{\lambda}]$, the canonical map 
     \[
       \colim_{\lambda \to \mu} \Map_{\catC_{\mu}}(c_{\mu}, c'_{\mu}) \to \Map_{\colim_{\Lambda} \catC_{\bullet}}(c, c')
     \]
     is an equivalence.
     Here, we write $c_{\mu}$ for the image of $c_{\lambda}$ under the transition functor $\catC_{\lambda} \to \catC_{\mu}$ and similarly for $c'_{\mu}$.
  \end{enumerate}
\end{recollection}

\begin{lemma}
  \label{lem:filtered-colimits-ff-and-ess-surj}
  Let $\Lambda$ be a filtered $\infty$-category, $\catC_{\bullet}, \catD_{\bullet} \from \Lambda \to \Cat_{\infty}$ two diagrams, and $f_{\!\bullet} \from \catC_{\bullet} \to \catD_{\bullet}$ a natural transformation.
  \begin{enumerate}[(1)]
    \item\label{lemitem:ess-surj} If for each $\lambda \in \Lambda$, the functor $f_{\!\lambda}$ is essentially surjective, so is $\colim_{\lambda \in \Lambda} f_{\!\lambda}$.
    \item\label{lemitem:ff} If for each $\lambda \in \Lambda$, the functor $f_{\!\lambda}$ is fully faithful, so is $\colim_{\lambda \in \Lambda} f_{\!\lambda}$.
  \end{enumerate}
\end{lemma}

\begin{proof}
  By \kerodon{068G}, the various functors $\{\catC_{\lambda} \to \colim_{\lambda \in \Lambda} \catC_{\lambda}\}_{\lambda \in \Lambda}$ are jointly surjective, and similarly for $\catD_{\bullet}$.
  Since for each $\lambda \in \Lambda$ the diagram
  \[
    \begin{tikzcd}
      \catC_{\lambda} \arrow[rr, "{f_{\!\lambda}}"] \arrow[d] && \catD_{\lambda} \arrow[d] \\
      \colim_{\lambda \in \Lambda} \catC_{\lambda} \arrow[rr, "{\colim_{\lambda \in \Lambda} f_{\!\lambda}}"'] && \colim_{\lambda \in \Lambda} \catD_{\lambda}
    \end{tikzcd}
  \]
  commutes, this shows~\ref{lemitem:ess-surj}.
  A proof of assertion~\ref{lemitem:ff} can be found in \cite[Prop. 2.1]{haine2025c}.
\end{proof}

\begin{lemma}
  \label{lem:limits-in-filtered-colimits-of-infinity-categories}
  Let $\catC_{\bullet} \from \Lambda \to \Cat_{\infty}, \lambda \mapsto \catC_{\lambda}$ be a filtered diagram of $\infty$-categories, and $K, \catE$ any other $\infty$-categories.
  For a map $\mu \to \lambda \in \Lambda$, let $f_{\!\mu \to \lambda} \from \catC_{\mu} \to \catC_{\lambda}$ denote the transition functor.
  \begin{enumerate}[(1)]
    \item \label{lemitem:colimits} Assume that $\catC_{\lambda}$ admits $K$-indexed colimits for each $\lambda \in \Lambda$ and that all transition functors preserve $K$-indexed colimits.
    \begin{enumerate}[(a)]
      \item\label{lemitem:existence-of-colimits} The $\infty$-category $\colim_{\Lambda}\!\catC_{\bullet}$ admits $K$-indexed colimits.
      \item \label{lemitem:inclusion-preserves-colimits} For each $\lambda \in \Lambda$ the canonical functor $i_{\lambda} \from \catC_{\lambda} \to \colim_{\Lambda}\!\catC_{\bullet}$ preserves $K$-indexed colimits.
      \item \label{lemitem:preservation-of-colimits} A functor $\colim_{\Lambda}\!\catC_{\bullet} \to \catE$ preserves $K$-indexed colimits if and only if the composition
      \[
        \catC_{\lambda} \to \colim_{\Lambda} \catC_{\bullet} \to \catE
      \]
      does for any $\lambda \in \Lambda$.
    \end{enumerate}
    \item \label{lemitem:finite-limits} Assume that $\catC_{\lambda}$ admits $K$-indexed limits for each $\lambda \in \Lambda$, and that all transition functors preserve $K$-indexed limits.
      Moreover, assume that $K$ is \emph{finite}.
    \begin{enumerate}[(a)]
        \item\label{lemitem:existence-of-finite-limits2} The $\infty$-category $\colim_{\Lambda}\! \catC_{\bullet}$ admits $K$-indexed limits.
        \item \label{lemitem:inclusion-preserves-finite-limits} For each $\lambda \in \Lambda$ the canonical functor $i_{\lambda} \from \catC_{\lambda} \to \colim_{\Lambda}\!\catC_{\bullet}$ preserves $K$-indexed limits.
        \item \label{lemitem:preservation-of-finite-limits} A functor $\colim_{\Lambda}\! \catC_{\bullet} \to \catE$ preserves $K$-indexed limits if and only if the composition
        \[
          \catC_{\lambda} \to \colim_{\Lambda}\! \catC_{\bullet} \to \catE
        \]
        does for any $\lambda \in \Lambda$.
    \end{enumerate}
  \end{enumerate}
\end{lemma}

\begin{proof}
  The case of colimits is a combination of \kerodon{0632} and \kerodon{0635}.
  For the rest of the proof, we assume that we are in situation \lref{lemitem:finite-limits}.
  Since $\Lambda$ is assumed to be filtered, any finite diagram in $\colim_{\Lambda}\! \catC_{\bullet}$ lifts to some finite stage $\lambda \in \Lambda$ by \kerodon{068G}.
  This shows that both \lref{lemitem:existence-of-finite-limits2} and \lref{lemitem:preservation-of-finite-limits} follow from \lref{lemitem:inclusion-preserves-finite-limits}.
  We now show that $i_{\lambda} \from \catC_{\lambda} \to \colim_{\Lambda}\! \catC_{\bullet}$ preserves $K$-indexed limits.
  Let $K \to \catC_{\lambda}, k \mapsto c_{\lambda, k}$ be a diagram and $c_{\Lambda} = \sqrbr{c_{\lambda'}} \in \colim_{\Lambda}\!\catC_{\bullet}$ any object.
  Since $\Lambda$ is filtered, we can find some $0 \in \Lambda$ and maps $\lambda \to 0 \ot \lambda'$ so that $c_{\Lambda} = \sqrbr{c_{0}}$ with $c_{0} = f_{\!\lambda' \to 0}(c_{\lambda'})$.
  Then, by \Cref{rec:colimits-of-infinity-categories}~\ref{recitem:maps-in-filtered-colimits}, we have
  \[
    \begin{split}
      \Map_{\colim_{\Lambda}\!\catC_{\bullet}}(c_{\Lambda}, i_{\lambda}(\limit_{k \in K} c_{\lambda, k})) &= \colim_{0 \to \mu} \Map_{\catC_{\mu}}(f_{\!0 \to \mu}(c_{0}), f_{\!\lambda \to 0 \to \mu}(\lim_{k \in K} c_{\lambda, k})) \\
                                                                  &= \colim_{0 \to \mu} \lim_{k \in K} \Map_{\catC_{\mu}}(f_{\!0 \to \mu}(c_0), f_{\!\lambda \to 0 \to \mu}(c_{\lambda, k})),
    \end{split}
  \]
  where the last equality holds by the assumption that each of the transition functors $f_{\!\lambda \to 0 \to  \mu}$ preserves $K$-indexed limits.
  The claim now follows from the fact that filtered colimits commute with finite limits in $\Ani$, see \kerodon{05XS}.
\end{proof}

\begin{lemma}
  \label{lem:filtered-colimits-and-slices}
  Let $\catC_{\bullet} \from \Lambda \to \Cat_{\infty}, \lambda \mapsto \catC_{\lambda}$ be a filtered diagram of $\infty$-categories with colimit $\catC_{\infty} \defeqq \colim_{\Lambda} \catC_{\bullet}$.
  Let $c_{\infty} \in \catC_{\infty}$ be obtained as the colimit of a compatible system of objects $c_{\lambda} \from \Delta^{0} \to \catC_{\lambda}$.
  The canonical map
  \[
    \colim_{\lambda \in \Lambda} (\catC_{\lambda})_{c_{\lambda}/} \to (\catC_{\infty})_{c_{\infty}/}
  \]
  is an equivalence of $\infty$-categories.
\end{lemma}

\begin{proof}
  Recall that, given an $\infty$-category $\catE$ and an object $e \in \catE$, the slice $\catE_{e/}$ is given by the pullback in $\Cat_{\infty}$
  \[
    \begin{tikzcd}
      \catE_{e/} \arrow[r] \arrow[d] \cartesian & \Fun(\Delta^{1}, \catE) \arrow[d, "{s}"] \\
      \Delta^{0} \arrow[r, "{e}"] & \catE,
    \end{tikzcd}
  \]
  where $s$ denotes the source functor.
  Thus, as $\Delta^{1}$ is compact in $\Cat_{\infty}$ by \Cref{rec:colimits-of-infinity-categories}~\ref{recitem:compact-categories}, the claim follows since filtered colimits commute with finite limits in $\Cat_{\infty}$.
\end{proof}


\section{Continuous local systems}
\label{sec:continuous-local-systems}
In this section, we study continuous local systems on profinite anima.

\begin{recollection}[continuous functors]
  \label{rec:currying}
  Let $\catC, \catD$ and $\catE$ be $\infty$-categories.
  \begin{enumerate}[(1)]
  \item\label{recitem:currying} 
    There is a canonical equivalence of $\infty$-categories
    \[
      \Fun(\catC \times \catD, \catE) \isomto \Fun(\catC, \Fun(\catD, \catE))
    \]
    given by $F \mapsto (c \mapsto F(c, \blank))$.
  \item \label{recitem:Fun}
  Pre- and postcomposition induce a functor
  \[
    \Fun \from \! \Cat_{\infty}^{\op} \times \Cat_{\infty} \to \Cat_{\infty}, (\catC, \catD) \mapsto \Fun(\catC, \catD).
  \]
  \item \label{recitem:Pro-Fun}
  Since $\Cat_{\infty}$ is cocomplete, so is $\Fun(\Cat_{\infty}, \Cat_{\infty})$.
  Therefore, the functor
  \[
    \Cat_{\infty}^{\op} \to \Fun(\Cat_{\infty}, \Cat_{\infty})
  \]
  corresponding under~\ref{recitem:currying} to the functor in~\ref{recitem:Fun} extends along filtered colimits to a functor
  \[
    \Pro(\Cat_{\infty})^{\op} = \Ind(\Cat_{\infty}^{\op}) \to \Fun(\Cat_{\infty}, \Cat_{\infty}).
  \]
  \item \label{recitem:continuous-functors}
  We write
  \[
    \Fun^{\cts} \from \Pro(\Cat_{\infty})^{\op} \times \Cat_{\infty} \to \Cat_{\infty}
  \]
  for the functor corresponding under~\ref{recitem:currying} to the functor constructed in~\ref{recitem:Pro-Fun}.
  Given a map $f \in \Pro(\Cat_{\infty})$, we write $f^{\ast} = (\blank \circ f)$ for the functor induced by precomposition.
  Given a functor $g$, we write $g_{\ast} = (g \circ \blank)$ for the functor induced by postcomposition.
  \item \label{recitem:evaluation}
  By construction, if $\catC = \prolimit_{\lambda \in \Lambda} \catC_{\lambda} \in \Pro(\Cat_{\infty})$, then
  \[
    \Fun^{\cts}(\catC, \catD) = \colim_{\lambda \in \Lambda^{\op}} \Fun(\catC_{\lambda}, \catD),
  \]
  where the transition functors are given by precomposition.
  We call $\Fun^{\cts}(\catC, \catD)$ the \emph{$\infty$-category of continuous functors from $\catC$ to $\catD$}.
\end{enumerate}
\end{recollection}

The inclusion $\Ani \subset \Cat_{\infty}$ induces an inclusion $\Pro(\Ani) \subset \Pro(\Cat_{\infty})$, allowing us to make the following definition:

\begin{definition}[continuous local systems]
  \label{def:continuous-local-systems}
  Let $B$ be any profinite anima, and $\catC$ an $\infty$-category.
  \begin{enumerate}[(1)]
    \item\label{defitem:continuous-C-valued-local-systems}
    We call $\Fun^{\cts}(B, \catC)$ the $\infty$-category of \emph{continuous $\catC$-valued local systems} on $B$.

    \item \label{defitem:continuous-local-systems}
    We simply refer to $\Fun^{\cts}(B, \Ani_{\pi})$ as the $\infty$-category of \emph{continuous local systems} on $B$.
  \end{enumerate}
\end{definition}

\subsection{(Co)Limits of continuous local systems}
\label{subsec:co-limits}

\begin{proposition}
  \label{lem:composition-preserves-finite-limits}
  Let $f \from \catC = \prolimit_{\lambda \in \Lambda} \catC_{\lambda} \to \prolimit_{\lambda \in \Lambda} \catC'_{\lambda} = \catC'$ be a map in $\Pro(\Cat_{\infty})$ and $g \from \catE \to \catE'$ a left-exact functor of $\infty$-categories with finite limits.
  Given any $\lambda \in \Lambda$, write $p_{\lambda} \from \catC \to \catC_{\lambda}$ and $p'_{\lambda} \from \catC' \to \catC'_{\lambda}$ for the projection maps.
  Let $K$ be an $\infty$-category.
  \begin{enumerate}[(1)]

    \item \label{lemitem:continuous-limits}
    Assume that $K$ is finite and that $\catE$ admits $K$-indexed limits.
    \begin{enumerate}[(a)]
      \item \label{lemitem:existence-of-finite-limits}
      The $\infty$-category $\Fun^{\cts}(\catC, \catE)$ admits $K$-indexed limits.
  
      \item \label{lemitem:preserves-limits}
      The functor $p_{\lambda}^{\ast} \from \Fun(\catC_{\lambda}, \catE) \to \Fun^{\cts}(\catC, \catE)$ preserves $K$-indexed limits.
  
      \item\label{lemitem:pre-preserves-limits}
      The functor $f^{\ast} \from \Fun^{\cts}(\catC', \catE) \to \Fun^{\cts}(\catC, \catE)$ preserves $K$-indexed limits.
  
      \item\label{lemitem:post-preserves-limits}
      If $g$ preserves $K$-indexed limits, so does $g_{\ast} \from \Fun^{\cts}(\catC, \catE) \to \Fun^{\cts}(\catC, \catE')$.
    \end{enumerate}

    \item \label{lemitem:continuous-colimits}
    Assume that $\cE$ admits $K$-indexed colimits.
    \begin{enumerate}[(a)]
      \item \label{lemitem:existence-of-finite-colimits}
      The $\infty$-category $\Fun^{\cts}(\catC, \catE)$ admits $K$-indexed colimits.
  
      \item \label{lemitem:preserves-colimits}
      The functor $p_{\lambda}^{\ast} \from \Fun(\catC_{\lambda}, \catE) \to \Fun^{\cts}(\catC, \catE)$ preserves $K$-indexed colimits.
  
      \item\label{lemitem:pre-preserves-colimits}
      The functor $f^{\ast} \from \Fun^{\cts}(\catC', \catE) \to \Fun^{\cts}(\catC, \catE)$ preserves $K$-indexed colimits.
  
      \item\label{lemitem:post-preserves-colimits}
      If $g$ preserves $K$-indexed colimits, so does $g_{\ast} \from \Fun^{\cts}(\catC, \catE) \to \Fun^{\cts}(\catC, \catE')$.
    \end{enumerate}
  \end{enumerate}
\end{proposition}

\begin{proof}
  Assertions~\ref{lemitem:continuous-limits}\ref{lemitem:existence-of-finite-limits} and~\ref{lemitem:preserves-limits} are special cases of \Cref{lem:limits-in-filtered-colimits-of-infinity-categories}~\ref{lemitem:finite-limits}\ref{lemitem:existence-of-finite-limits2} and \lref{lemitem:inclusion-preserves-finite-limits}, respectively.
  Next, we verify~\ref{lemitem:continuous-limits}\ref{lemitem:pre-preserves-limits} and~\ref{lemitem:post-preserves-limits}.
  By construction, the diagrams
  \[
    \begin{tikzcd}
      \Fun(\catC'_{\lambda}, \catE) \arrow[d, "{(p'_{\lambda})^{\ast}}"'] \arrow[r, "{f^{\ast}}"] & \Fun(\catC_{\lambda}, \catE) \arrow[d, "{p_{\lambda}^{\ast}}"] & \Fun(\catC_{\lambda}, \catE) \arrow[d, "{p_{\lambda}^{\ast}}"'] \arrow[r, "{g_{\ast}}"] & \Fun(\catC_{\lambda}, \catE') \arrow[d, "{p_{\lambda}^{\ast}}"] \\
      \Fun^{\cts}(\catC', \catE) \arrow[r, "{f^{\ast}}"'] & \Fun^{\cts}(\catC, \catE) &  \Fun^{\cts}(\catC, \catE) \arrow[r, "{g_{\ast}}"'] & \Fun^{\cts}(\catC, \catE')
    \end{tikzcd}
  \]
  commute for each $\lambda \in \Lambda$.
  Therefore, by \Cref{lem:limits-in-filtered-colimits-of-infinity-categories}~\ref{lemitem:finite-limits}\lref{lemitem:inclusion-preserves-finite-limits} and \lref{lemitem:preservation-of-finite-limits}, we are reduced to the analogous statements for ordinary functor categories.
  Here, the result is obvious because (co)limits in functor categories are computed pointwise, see (the dual of) \kerodon{02X9}.
  The argument for~\ref{lemitem:continuous-colimits} is analogous.
\end{proof}

\begin{lemma}
  \label{lem:continuous-adjunctions}
  Let $L \from \catE \to \catE'$ be a functor and $\catC \in \Pro(\Cat_{\infty})$.
  \begin{enumerate}[(1)]
  \item\label{lemitem:cts-adjunction} If $R$ is right-adjoint to $L$, then the pair of functors
    \[
      \begin{tikzcd}
        \Fun^{\cts}(\catC, \catE) \arrow[r, "{L_{\ast}}", shift left = 0.3em] & \Fun^{\cts}(\catC, \catE') \arrow[l, "{R_{\ast}}", shift left = 0.3em]
      \end{tikzcd}
    \]
    induced by postcomposition are again adjoint to each other.
  \item \label{lemitem:cts-fully-faithful}
    If $L$ is fully faithful, so is $L_{\ast} \from \Fun^{\cts}(\catC, \catE) \to \Fun^{\cts}(\catC, \catE')$.
  \end{enumerate}
\end{lemma}

\begin{proof}
  Choose a cofiltered diagram $\catC_{\bullet} \from \Lambda \to \Cat_{\infty}, \lambda \mapsto \catC_{\lambda}$ representing $\catC$ and a unit $u \from \id_{\catE} \to R \circ L$ exhibiting $R$ as right adjoint to $L$.
  Write $i_{\lambda} \from \Fun(\catC_{\lambda}, \catE) \to \Fun^{\cts}(\catC, \catE)$ for the canonical inclusion functor.
  Assertion~\ref{lemitem:cts-adjunction} holds for the ordinary functor categories on each finite stage $\lambda \in \Lambda$.
  More precisely, $u$ induces a unit $u^{\lambda} \from \id_{\Fun(\catC_{\lambda}, \catE)} \to R_{\ast} \circ L_{\ast}$ with components $u^{\lambda}_F = (u_{F(c_{\lambda})})_{c_{\lambda} \in \catC_{\lambda}}$ exhibiting $R_{\ast} \from \Fun(\catC_{\lambda}, \catE') \to \Fun(\catC_{\lambda}, \catE)$ as right adjoint of the corresponding $L_{\ast}$.
  Moreover, since the transition functors in $\{\Fun(\catC_{\lambda}, \catE)\}_{\lambda \in \Lambda}$ are given via \emph{pre}composition, the various functors induced by \emph{post}composition with $L$ assemble into a natural transformation $L_{\bullet} \from \Fun(\catC_{\bullet}, \catE) \to \Fun(\catC_{\bullet}, \catE') $, and similarly for $R$.
  By definition, $L_{\ast}$ (resp. $R_{\ast}$) is the colimit of $L_{\bullet}$ (resp. $R_{\bullet}$), and
  \[
    (i_{\lambda}u^{\lambda})_{\lambda \in \Lambda} \in \lim_{\lambda \in \Lambda}\Fun(\Fun(\catC_{\lambda}, \catE), \Fun^{\cts}(\catC_{\lambda}, \catE)) = \Fun(\Fun^{\cts}(\catC, \catE), \Fun^{\cts}(\catC, \catE))
  \]
  defines a unit exhibiting that $R_{\ast} \from \Fun^{\cts}(\catC, \catE') \to \Fun^{\cts}(\catC, \catE)$ is right adjoint to $L_{\ast}$ as claimed.  
  Assertion~\ref{lemitem:cts-fully-faithful} follows from the fact that the analogous statement is true for ordinary functor categories, together with \Cref{lem:filtered-colimits-ff-and-ess-surj}~\ref{lemitem:ff}.
\end{proof}

\subsection{Truncations of continuous local systems}
\label{subsec:truncations-of-continuous-local-systems}

\begin{notation}
  \label{not:truncated-objects}
  Let $\catC$ be an $\infty$-category and $k \geq -1$.
  We write $\catC_{\leq k}$ for the full subcategory of $\catC$ spanned by the $k$-truncated objects.
\end{notation}

\begin{notation}
  \label{not:evaluation}
  Let $\catC \in \Pro(\Cat_{\infty})$ and $\catE$ an $\infty$-category.
  \begin{enumerate}[(1)]
    \item\label{notitem:object} An \emph{object} of $\catC$ is a map $c \from \Delta^0 \to \catC$ in $\Pro(\Cat_{\infty})$.
      We usually simply write $c \in \catC$ in this case.
    \item \label{notitem:evaluation} Let $c \in \catC$ be an object.
      We write $\ev_{\! c} \defeqq c^{\ast} \from \Fun^{\cts}(\catC, \catE) \to \Fun^{\cts}(\Delta^0, \catE) = \catE$ and refer to $\ev_{\! c}$ as the \emph{evaluation at $c$}.
      \item \label{notitem:constant-diagram} Let $e \in \catE$ be an object.
      We write $\const_e \in \Fun^{\cts}(\catC, \catE)$ for the image of $e$ under the functor $\catE = \Fun^{\cts}(\Delta^0, \catE) \to \Fun^{\cts}(\catC, \catE)$ induced by precomposition with the unique map $\catC \to \Delta^{0} \in \Pro(\Cat_{\infty})$ and refer to $\const_e$ as the \emph{constant diagram with value $e$}.
  \end{enumerate}
\end{notation}

\begin{lemma}
  \label{lem:continuous-functors-into-n-category}
  Let $B \in \Pro(\Ani_{\pi})$, $n \geq 0$ and $\catE$ an $n$-category.
  Precomposition with the $n$-truncation map $B \to \tau_{\leq n} B$ induces an equivalence $\Fun^{\cts}(\tau_{\leq n}B, \catE) \to \Fun^{\cts}(B, \catE)$.
\end{lemma}

\begin{proof}
  Write $B = \prolimit_{\lambda \in \Lambda} B_{\lambda}$ as a cofiltered limit of $\pi$-finite anima.
  Then $B \to \tau_{\leq n}B$ is the cofiltered limit of the $n$-truncation maps $B_{\lambda} \to \tau_{\leq n}B_{\lambda}$.
  We may thus assume that $B \in \Ani_{\pi}$ and thus $\Fun^{\cts}(B, \catE) = \Fun(B, \catE)$.
  Since $\catE$ is an $n$-category, it is an $n$-truncated object of $\Cat_{\infty}$.
  The claim now follows from the fact that the inclusion $\Ani \subset \Cat_{\infty}$ preserves $n$-truncations.
\end{proof}

\begin{proposition}
  \label{prop:continuous-truncations}
  Let $\catC \in \Pro(\Cat_{\infty})$, $B \in \Pro(\Ani_{\pi})$, $\catE$ an $\infty$-category and $k \geq -1$ and integer.
  \begin{enumerate}[(1)]
  \item\label{propitem:cts-k-truncated-objects} 
  The inclusion $\catE_{\leq k} \subset \catE$ induces a functor
  \[
    \Fun^{\cts}(\catC, \catE_{\leq k}) \to \Fun^{\cts}(\catC, \catE)
  \]
  that is fully faithful with essential image contained in $\Fun^{\cts}(\catC, \catE)_{\leq k}$.

  \item \label{propitem:ess-surjectivity} We have $\Fun^{\cts}(B, \catE_{\leq k}) = \Fun^{\cts}(B, \catE)_{\leq k}$.
  \item \label{propitem:cts-k-truncated}
  If the inclusion $\catE_{\leq k} \subset \catE$ admits a left-adjoint $\tau_{\leq k} \from \catE \to \catE_{\leq k}$, then $(\tau_{\leq k} \circ \blank)$ is left-adjoint to the inclusion $\Fun^{\cts}(B, \catE)_{\leq k} \subset \Fun^{\cts}(B, \catE)$.
  \end{enumerate}
\end{proposition}

\begin{proof}
  By \Cref{lem:continuous-adjunctions}~\ref{lemitem:cts-fully-faithful}, the only thing left to show in~\ref{propitem:cts-k-truncated-objects} is that $\Fun^{\cts}(\catC, \catE_{\leq k}) \subset \Fun^{\cts}(\catC, \catE)_{\leq k}$.
  This follows from \Cref{lem:composition-preserves-finite-limits}~\ref{lemitem:continuous-limits}\ref{lemitem:post-preserves-limits} combined with \HTT{}{5.5.6.5} and \HTT{}{5.5.6.16}.
  To verify~\ref{propitem:ess-surjectivity}, pick a level-representation $B = \prolimit_{\lambda \in \Lambda} B_{\lambda}$ with the property that each of the canonical projections $p_{\lambda} \from B \to B_{\lambda}$ is surjective on connected components as in \Cref{lem:level-representation-profinite-anima}.
  Let $L \in \Fun^{\cts}(B, \catE)_{\leq k}$.
  Then there exists some $\lambda \in \Lambda$ and $L_{\lambda} \in \Fun(B_{\lambda}, \catE)$ such that $L = L_{\lambda} \circ p_{\lambda}$.
  It therefore suffices to conclude that $L_{\lambda}$ takes values in $\catE_{\leq k}$.
  Let $b_{\lambda} \in B_{\lambda}$ and choose some $b \in B$ with the property that $p_{\lambda}(b)$ lies in the same connected component of $B_{\lambda}$ as $b_{\lambda}$.
  Then $L_{\lambda}(b_{\lambda}) \cong \ev_b(L) \in \catE$.
  Therefore, in order to show that $L_{\lambda}$ takes values in $\catE_{\leq k} \subset \catE$, it suffices to show that $\ev_b(L) \in \catE_{\leq k}$ for each $b \in B$.
  This follows from \Cref{lem:composition-preserves-finite-limits}~\ref{lemitem:continuous-limits}\ref{lemitem:pre-preserves-limits} combined with \HTT{}{5.5.6.5} and \HTT{}{5.5.6.16}.
  Finally,~\ref{propitem:cts-k-truncated} follows from~\ref{propitem:ess-surjectivity} combined with \Cref{lem:continuous-adjunctions}~\ref{lemitem:cts-adjunction}.
\end{proof}

\subsection{Relation with classical local systems}
\label{subsec:relation-with-classical-local-systems}

In this section, we relate continuous local systems of $\pi$-finite anima with continuous versions of classical local systems.

\begin{lemma}
  \label{lem:slice-of-cts-functor-category}
  Let $\catC \in \Pro(\Cat_{\infty})$, $\catE$ an $\infty$-category, $e \in \catE$ an object, and write $\const_e \in \Fun(\catC, \catE)$ for the constant diagram as in \Cref{not:evaluation} \lref{notitem:constant-diagram}.
  Postcomposition with the forgetful functor $V \from \catE_{e/} \to \catE$ induces an equivalence of $\infty$-categories $\Fun^{\cts}(\catC, \catE_{e/}) \to \Fun^{\cts}(\catC, \catE)_{\const_{e}\!/}$.
\end{lemma}

\begin{proof}
  Since $\const_{\id_e}$ defines an initial object of $\Fun^{\cts}(\catC, \catE_{e/})$ and $V(\const_{\id_e}) = \const_e$, we get the induced functor.
  To see that it is an equivalence, note that \Cref{lem:filtered-colimits-and-slices} reduces us to the case where $\catC \in \Cat_{\infty}$.
  Now recall that both diagrams 
  \[
    \begin{tikzcd}
      \catE_{e/} \arrow[r] \arrow[d] \cartesian & \Fun(\Delta^1, \catE) \arrow[d, "{s}"] & \Fun(\catC, \catE)_{\const_e\!/} \arrow[r] \arrow[d] \cartesian & \Fun(\Delta^1, \Fun(\catC, \catE)) \arrow[d, "s"]  \\
      \Delta^0 \arrow[r, "{e}"'] & \catE & \Delta^0 \arrow[r, "{\const_{e}}"'] \arrow[r] & \Fun(\catC, \catE)
    \end{tikzcd}
  \]
  are cartesian and similarly for $\Fun(\catC, \catE)_{\const_e\!/}$.
  The claim now follows from applying and carefully unwinding the chain of equivalences $\Fun(\Delta^1, \Fun(\catC, \catE)) \simeq \Fun(\Delta^1 \times \catC, \catE) \simeq \Fun(\catC, \Fun(\Delta^1, \catE))$.
\end{proof}

\begin{notation}
  Let $\catC$ be an $\infty$-category with terminal object $\ast \in \catC$.
  We write $\catC_{\ast} \defeqq \catC_{\ast/}$ for the $\infty$-category of pointed objects in $\catC$.
\end{notation}

\begin{proposition}
  \label{prop:pointed--and-group-objects}
  Let $\catC \in \Pro(\Cat_{\infty})$ and $\catE$ an $\infty$-category.
  \begin{enumerate}[(1)]
    \item\label{propitem:pointed-objects} 
    If $\ast \in \catE$ is a terminal object, the forgetful functor $V \from \catE_{\ast} \to \catE$ determines an equivalence of $\infty$-categories $V_{\!\ast} \from \Fun^{\cts}(\catC, \catE_{\ast}) \isomto \Fun^{\cts}(\catC, \catE)_{\ast}$.
    \item \label{propitem:groupoid-objects}
    If $\catE$ admits pullbacks, we have a canonical equivalence of $\infty$-categories
    \[
      \Fun^{\cts}(\catC, \Grpd(\catE)) \isomto \Grpd(\Fun^{\cts}(\catC, \catE)).
    \]
    \item \label{propitem:group-objects}
    If $\catE$ admits finite products, we have a canonical equivalence of $\infty$-categories
    \[
      \Fun^{\cts}(\catC, \Grp(\catE)) \isomto \Grp(\Fun^{\cts}(\catC, \catE)). 
    \]
  \end{enumerate}
  
\end{proposition}

\begin{proof}
  Since $\ast \in \catE_{\ast}$ is initial, so is $\const_{\ast} \in \Fun^{\cts}(\catC, \catE_{\ast})$.
  Moreover, $V(\const_{\ast})$ is a terminal object of $\Fun^{\cts}(\catC, \catE)$.
  Assertion \lref{propitem:pointed-objects} is thus an immediate consequence of \Cref{lem:slice-of-cts-functor-category}.
  Assertions \lref{propitem:groupoid-objects} and \lref{propitem:group-objects} are straightforward consequences of \cite[\SAGthm{}{A.8.2.5}]{SAG}.
\end{proof}

In summary, we have the following relation between $\Ani_{\pi}$-valued and classical local systems:

\begin{theorem}
  \label{thm:relation-to-classical-local-systems}
  Let $B$ be a profinite anima.
  \begin{enumerate}[(1)]
    \item\label{thmitem:n-truncations}
    Postcomposition with the inclusion $\Ani_{\pi, \leq n} \subset \Ani_{\pi}$ induces a fully faithful embedding
    \[
      \Fun^{\cts}(B, \Ani_{\pi, \leq n}) \subset \Fun^{\cts}(B, \Ani_{\pi})
    \]
    with essential image given by the $n$-truncated objects.
    Moreover, this inclusion admits a left-adjoint obtained by postcomposition with the usual $n$-truncation functor $\tau_{\leq n} \from \Ani_{\pi} \to \Ani_{\pi, \leq n}$.
    \item \label{thmitem:0-truncation}
    The inclusion $\Set_{\fin} = \Ani_{\pi, \leq 0} \subset \Ani_{\pi}$ induces an equivalence
    \[
      \Fun^{\cts}(B, \Set_{\fin}) \isomto \Fun^{\cts}(B, \Ani_{\pi})_{\leq 0}.
    \]
    \item \label{thmitem:pointed-objects}
    There is a canonical equivalence of $\infty$-categories
    \[
      \Fun^{\cts}(B, \Ani_{\pi, \ast}) \isomto \Fun^{\cts}(B, \Ani_{\pi})_{\ast}.
    \]
    \item \label{thmitem:groups}
    We have a canonical equivalence of $1$-categories
    \[
      \Fun^{\cts}(B, \Grp_{\fin}) \isomto \Grp(\Fun^{\cts}(B, \Set_{\fin})).
    \]
    \item \label{thmitem:abelian-groups}
    We have a canonical equivalence of $1$-categories
    \[
      \Fun^{\cts}(B, \Ab_{\fin}) \isomto \Ab(\Fun^{\cts}(B, \Set_{\fin})).
    \]
  \end{enumerate}
\end{theorem}

\begin{proof}
  Assertion \lref{thmitem:n-truncations} is a special case of \Cref{prop:continuous-truncations}, \lref{thmitem:0-truncation} is a special case of \lref{thmitem:n-truncations}, \lref{thmitem:pointed-objects} is a special case of \Cref{prop:pointed--and-group-objects} \lref{propitem:pointed-objects}, \lref{thmitem:groups} is a special case of \Cref{prop:pointed--and-group-objects} \lref{propitem:pointed-objects}, and, since $\Ab = \Grp(\Grp)$, \lref{thmitem:abelian-groups} follows by applying \Cref{prop:pointed--and-group-objects} \lref{propitem:pointed-objects} twice.
\end{proof}


\section{Continuous (un-)straightening}
\label{sec:continous-un-straightening}
The goal of this section is to prove \Cref{thm:continuous-un-straightening}, a continuous analogue of the following.

\begin{recollection}
  \label{rec:un-straightening-for-anima}
  Let $B \in \Ani$.
    There is an adjoint equivalence of $\infty$-categories
    \[
      \begin{tikzcd}
        \Unstr^B \from \Fun(B, \Ani) \arrow[r, shift left = 0.2em] & \Ani_{/B} \arrow[l, shift left = 0.2em] \from \Str_B
      \end{tikzcd}
    \]
    called the \emph{straightening-unstraightening equivalence}.
    Here, $\int^{B} \defeqq \colim_B$ denotes the colimit functor and $\Str_{B}(f \from E \to B)$ is the local system of fibres $b \mapsto \fib_b(f)$.
\end{recollection}

Lurie shows in \cite[\SAGapp{E}]{SAG} that the colimit functor
\[
  \colim_B \from \Fun(B, \Pro(\Ani_{\pi})) \to \Pro(\Ani_{\pi})_{/B}
\]
is an equivalence if $B$ is $\pi$-finite, but not for general profinite anima (see \cites[\SAGthm{}{E.6.0.7}, \SAGthm{}{E.6.0.9}]{SAG}).
In order to get an equivalence for profinite anima, one has to replace $\Fun(B, \Pro(\Ani_{\pi}))$ with the $\infty$-category $\Pro(\Fun^{\cts}(B, \Ani_{\pi}))$.

In fact, we will prove a continuous relative pro-$\Sigma$ variant of \Cref{rec:un-straightening-for-anima}.
To this end, we introduce the following notion.

\begin{definition}
  \label{def:Sigma-finite-maps-finite-case}
  Let $\Sigma$ denote a full class of finite groups and $f \from E \to B$ a map of $\pi$-finite anima.
  \begin{enumerate}
    \item\label{defitem:Sigma-finite-map-finite-case}
    The map $f$ is said to be \emph{$\Sigma$-finite} if for every point $b \in B$, we have that $\fib_b(f) \in \Ani_{\Sigma}$.
    \item \label{defitem:category-of-Sigma-finite-maps}
    We write $\Ani_{\Sigma} \rel B$ for the full subcategory of $\Ani_{\pi / B}$ spanned by the $\Sigma$-finite maps.
  \end{enumerate}
\end{definition}

\begin{lemma}[{$\Sigma$-finite straightening-unstraightening}]
  \label{lem:straightening-unstraigthening}
  Let $\Sigma$ denote a nonempty full class of finite groups.
  \begin{enumerate}[(1)]
  \item\label{lemitem:2-out-of-3} If $B \in \Ani_{\Sigma}$, we have $\Ani_{\Sigma/B} = \Ani_{\Sigma}\rel B  \subset \Ani_{\pi/B} $.
    In particular $ \Ani_{\pi} \rel B = \Ani_{\pi/B}$.
    \item \label{lemitem:pi-finite-(un)-straightening}
    The straightening-unstraightening equivalence of \Cref{rec:un-straightening-for-anima}restricts to an adjoint equivalence  of $\infty$-categories
    \[
      \begin{tikzcd}
        \Fun(B, \Ani_{\Sigma}) \arrow[rr, shift left = 0.2em, "{\Unstr_{\!B}}"] & & \Ani_{\Sigma} \rel B . \arrow[ll, shift left = 0.2em, "{\Str_{B}}"]
      \end{tikzcd}
    \]   
    \item \label{lemitem:functoriality}
    For any map of $\pi$-finite anima $g \from B' \to B$, the diagram
    \[
      \begin{tikzcd}
        \Fun(B, \Ani_{\Sigma}) \arrow[d, "{\blank \circ g}"'] & \Ani_{\Sigma} \rel B \arrow[l, "{\Str_{B}}"'] \arrow[d, "{\blank \times_B B'}"] \\
        \Fun(B', \Ani_{\Sigma}) & \Ani_{\Sigma} \rel B' . \arrow[l, "{\Str_{B'}}"]
      \end{tikzcd}
    \]
    commutes.
  \end{enumerate}
\end{lemma}

\begin{proof}
  The inclusion $\Ani_{\Sigma/B} \subset \Ani_{\Sigma} \rel B$ is an immediate consequence of the fact that $\Ani_{\Sigma} \subset \Ani$ is closed under finite limits.
  If $E \in \Ani_{\Sigma}\rel B$, both $B$ and all fibres $E_b = E \times_B b$ are $\Sigma$-finite.
  By inspecting the induced long exact sequences on homotopy groups for varying $b \in B$ and $e \in E_b$, one readily checks that $E \in \Ani_{\Sigma}$.
  This shows~\ref{lemitem:2-out-of-3}.
  Assertion~\ref{lemitem:pi-finite-(un)-straightening} is an immediate consequence of~\ref{lemitem:2-out-of-3}, and~\ref{lemitem:functoriality} follows from \Cref{cor:fibrewise-criterion-pullbacks}.
\end{proof}

Next, we want to extend \Cref{def:Sigma-finite-maps-finite-case} to maps of profinite anima.
The following definition and lemma are essentially a variant of \cite[Lemma 4.13]{clausen2025}.

\begin{definition}[$\Sigma$-(pro-)finite maps]
  \label{def:Sigma-finite-maps}
  Let $\Sigma$ be a full class of finite groups and $f \from E \to B$ a map of profinite anima.
  \begin{enumerate}[(1)]
  \item\label{defitem:Sigma-finite-map} 
  The map $f \from E \to B$ is called \emph{$\Sigma$-finite} if there exists a cartesian square
  \[
    \begin{tikzcd}
      E \arrow[r] \arrow[d, "{f}"'] \cartesian & E_{\pi} \arrow[d, "{f_{\!\Sigma}}"] \\
      B \arrow[r] & B_{\pi}
    \end{tikzcd}
  \]
  in $\Pro(\Ani_{\pi})$ with $\pi$-finite anima $E_{\pi}, B_{\pi}$ and $f_{\!\Sigma} \in \Ani_{\Sigma} \rel B_{\pi}$.
  We write $\Ani_{\Sigma} \rel B \subset \Pro(\Ani_{\pi})_{/B}$ for the full subcategory spanned by the $\Sigma$-finite maps.

  \item \label{defitem:Sigma-profinite-map}
  The map $f \from E \to B$ is called \emph{$\Sigma$-profinite} if for any point $b \in B$, we have that $\fib_b(f) \in \Pro(\Ani_{\Sigma})$.
  We write $\Pro(\Ani_{\Sigma}) \rel B \subset \Pro(\Ani_{\pi})_{/B}$ for the full subcategory spanned by the $\Sigma$-profinite maps.
  \end{enumerate}
\end{definition}

\begin{proposition}
  \label{lem:Sigma-finite-maps}
  Let $B = \prolimit_{\lambda \in \Lambda} B_{\lambda}$ be a profinite anima.
  \begin{enumerate}[(1)]
    \item \label{lemitem:embedding-to-profinite-anima-over-B}
    The various functors $\Ani_{\Sigma} \rel B_{\lambda} \to \Pro(\Ani_{\pi})_{/B}$ given by $E_{\lambda} \mapsto B \times_{B_{\lambda}} E_{\lambda}$ assemble into a fully faithful functor
    \[
       \colim_{\lambda \in \Lambda^{\op}} (\Ani_{\Sigma}\rel B_{\lambda}) \to \Pro(\Ani_{\pi})_{/B}
    \]
    with essential image $\Ani_{\Sigma}\rel B$.
    \item \label{lemitem:Sigma-rel-Sigma}
    If $B \in \Pro(\Ani_{\Sigma})$, then $\Pro(\Ani_{\Sigma}) \rel B = \Pro(\Ani_{\Sigma})_{/B}$.
    \item \label{lemitem:equivalence-B-pro-Sigma}
    Extending along cofiltered limits yields an inclusion of $\infty$-categories
    \[
      \Pro(\Ani_{\Sigma}\rel B) \subset \Pro(\Ani_{\Sigma}) \rel B.
    \]
    that is an equivalence if $B \in \Pro(\Ani_{\Sigma})$.
  \end{enumerate}
\end{proposition}

\begin{proof}
  Since for any map $\mu \to \lambda \in \Lambda$, the diagram
  \[
    \begin{tikzcd}
      \Ani_{\Sigma} \rel B_{\lambda} \arrow[rd, "{\blank \times_{B_{\lambda}} B}"'] \arrow[rr, "{\blank \times_{B_{\lambda}} B_{\mu}}"] && \Ani_{\Sigma} \rel B_{\mu} \arrow[ld, "{\blank \times_{B_{\mu}} B}"] \\
      & \Pro(\Ani_{\pi})_{/B}
    \end{tikzcd}
  \]
  commutes, we get the functor in~\ref{lemitem:embedding-to-profinite-anima-over-B}.
  To see fully faithfulness, let $[E_{0} \to B_{0}]$ and $[E'_{0} \to B_{0}]$ be two objects in $\Ani_{\Sigma}\rel B$ that we may assume to be represented at some common level $0$ since $\Lambda$ is cofiltered.
  Then, by \Cref{rec:colimits-of-infinity-categories}~\ref{recitem:maps-in-filtered-colimits}, we have 
  \[
    \begin{split}
      \Map([E'_{0}], [E_{0}]) &= \colim_{0 \to \lambda \in \Lambda^{\op}} \Map_{B_{\lambda}}(E'_{0} \times_{B_{0}} B_{\lambda}, E_{0} \times_{B_{0}} B_{\lambda}) \\
                             &= \colim_{0 \to \lambda \in \Lambda^{\op}} \Map_{B_{0}}(E'_{0} \times_{B_{0}} B_{\lambda}, E_{0}) \\
                             &= \Map_{\Pro(\Ani_{\pi})_{/B_{0}}}(\prolimit_{\lambda \to 0 \in \Lambda}(E'_{0} \times_{B_{0}} B_{\lambda}), E_{0}) \\
                             &= \Map_{\Pro(\Ani_{\pi})_{/B}}(\prolimit_{\lambda \to 0 \in \Lambda}(E'_{0} \times_{B_{0}} B_{\lambda}), E_{0} \times_{B_{0}} B)
    \end{split}
  \]
  which, since $\prolimit_{\lambda \to 0 \in \Lambda}(E'_{0} \times_{B_{0}} B_{\lambda}) = E'_{0} \times_{B_{0}} B$ by the dual of \Cref{rec:colimits-of-infinity-categories}~\ref{recitem:colimit-cofinal}, shows fully faithfulness.
  Clearly, the image is contained in $\Ani_{\Sigma}\rel B$.
  To see equality, let $f \from E \to B$ be $\Sigma$-finite and choose a cartesian diagram
  \[
    \begin{tikzcd}
      E \arrow[r] \arrow[d, "{f}"'] \cartesian & E_{\pi} \arrow[d, "{f_{\!\Sigma}}"] \\
      B \arrow[r] & B_{\pi}
    \end{tikzcd}
  \]
  exhibiting the $\Sigma$-finiteness of $f$.
  Since $B_{\pi}\in \Ani_{\pi}$, the map $B \to B_{\pi}$ factors through one the $B_{\lambda}$'s.
  Then $E_{\lambda} \defeqq E_{\pi} \times_{f_{\!\Sigma}, B_{\pi}} B_{\lambda} \to B_{\lambda} \in \Ani_{\Sigma}\rel B_{\lambda}$ pulls back to $f$.
  This completes the proof of~\ref{lemitem:embedding-to-profinite-anima-over-B}.
  Next, we verify~\ref{lemitem:Sigma-rel-Sigma}.
  That $\Pro(\Ani_{\Sigma})_{/B} \subset \Pro(\Ani_{\Sigma}) \rel B$ follows from the fact that $\Pro(\Ani_{\Sigma}) \subset \Pro(\Ani_{\pi})$ is closed under pullbacks.
  To see equality, let $f \from E \to B$ be a map of profinite anima in $\Pro(\Ani_{\Sigma})\rel B$.
  We need to show that $E \in \Pro(\Ani_{\Sigma})$.
  To this end, let $b \in B$ be any point and consider the corresponding map of fibre sequences
  \[
    \begin{tikzcd}
      F_{b} \arrow[r] \arrow[d] & E \arrow[r, "f"] \arrow[d] & B \arrow[d, "{\cong}"] \\
      F_{b}(\Sigma) \arrow[r] & E\Spfcompl \arrow[r, "{f\Spfcompl}"'] & B\Spfcompl.
    \end{tikzcd}
  \]
  Since $B \in \Pro(\Ani_{\Sigma})$, the right vertical map is an equivalence.
  Moreover, we have
  \[
    F_b \isomto (F_b)\Spfcompl = (b \times_B E)\Spfcompl \isomto b \times_{B\Spfcompl} E\Spfcompl = F_{b}(\Sigma),
  \]
  where the first isomorphism comes from the assumption that $f \in \Pro(\Ani_{\Sigma})\rel B$, and the second isomorphism follows from \Cref{rec:Spfcompl-preserves-cofiltered-limits}~\ref{recitem:pfcompl-is-locally-cartesian-and-preserves-limits} together with the assumption that $B \in \Pro(\Ani_{\Sigma})$.
  Since this holds for arbitrary $b \in B$, we conclude that $E \to E\Spfcompl$ is an equivalence by \Cref{thm:profinite-whitehead}.
  Next, we verify~\ref{lemitem:equivalence-B-pro-Sigma}.
  That $\Ani_{\Sigma} \rel B \subset \Pro(\Ani_{\Sigma})\rel B$ is an immediate consequence of \Cref{cor:fibrewise-criterion-pullbacks}.
  To see that $\Pro(\Ani_{\Sigma}\rel B) \to \Pro(\Ani_{\Sigma})\rel B$ is fully faithful, it is enough to show that every object $E \in \Ani_{\Sigma} \rel B$ is cocompact in $\Pro(\Ani_{\pi})_{/B}$.
  To this end, simply note that for $E = E_{\lambda} \times_{B_{\lambda}} B$, we have that
  \[
    \Map_B(\blank, E) = \Map_{B_{\lambda}}(\blank, E_{\lambda}) = \Map(\blank, E_{\lambda}) \times_{\Map(\blank, B_{\lambda})} \ast
  \]
  takes cofiltered limits to filtered colimits since $E_{\lambda}, B_{\lambda} \in \Ani_{\pi}$ are cocompact in $\Pro(\Ani_{\pi})$ and since colimits are universal in $\Ani$.
  To finish the proof of~\ref{lemitem:equivalence-B-pro-Sigma}, the only thing left to show is that any map of profinite anima $f \from E \to B$ in $\Pro(\Ani_{\Sigma})\rel B$ can be written as a cofiltered limit of maps $f_{\!(\lambda)} \in \Ani_{\Sigma}\rel B$ if $B \in \Pro(\Ani_{\Sigma})$.
  With $B$, also $E \in \Pro(\Ani_{\Sigma})$ by~\ref{lemitem:Sigma-rel-Sigma}.
  We can thus write $f = \prolimit_{\lambda \in \Lambda} f_{\!\lambda}$ with $f_{\!\lambda} \from E_{\lambda} \to B_{\lambda}$ a map in $\Ani_{\Sigma/B_{\lambda}} = \Ani_{\Sigma} \rel B_{\lambda}$.
  Since $f = \limit_{\lambda} (f_{\!\lambda} \times_{B_{\lambda}} B)$, we conclude.
\end{proof}

\begin{corollary}
  \label{cor:pi-finite-maps-stable-under-composition}
  Let $f \from E' \to E$ and $g \from E \to B$ be maps of profinite anima.
  If $g \in \Ani_{\pi}\rel B$, then $f \in \Ani_{\pi}\rel E$, if and only if $(g \circ f) \in \Ani_{\pi} \rel B$.
\end{corollary}

\begin{proof}
  Write $E' \xrightarrow{f} E \xrightarrow{g} B$ as a cofiltered limit of maps $E'_{\lambda} \xrightarrow{f_{\!\lambda}} E_{\lambda} \xrightarrow{g_{\lambda}} B_{\lambda}$ over some cofiltered category $\lambda \in \Lambda$ and choose a pullback square
  \[
    \begin{tikzcd}
      E \arrow[r] \arrow[d, "g"'] \cartesian & E_{\pi} \arrow[d, "{g_{\pi}}"] \\
      B \arrow[r] & B_{\pi}  
    \end{tikzcd}
  \]
  exhibiting the $\pi$-finiteness of $g$.
  Since $B_{\pi}$ is $\pi$-finite, there exists some $\lambda \in \Lambda$ such that $B \to B_{\pi}$ factors through $B \to B_\lambda$.
  Similarly, since $E_{\pi}$ is $\pi$-finite and $\Lambda$ is cofiltered, we can find some $\mu \to \lambda \in \Lambda$ such that $E \to E_{\pi} \times_{B_{\pi}} B_{\lambda}$ factors through $E_\mu$.
  Replacing $\Lambda$ with $\Lambda_{/\mu}$, we may hence assume that for any $\lambda \in \Lambda$ we have a commutative diagram
  \[
    \begin{tikzcd}
      E' \arrow[d, "f"'] \arrow[r] & E'_{\lambda} \arrow[d, "f_{{\!\lambda}}"] \\
      E \arrow[r] \arrow[d, "g"'] & E_{\lambda} \arrow[r] \arrow[d, "{g_{\lambda}}"] & E_{\pi} \arrow[d, "{g_{\pi}}"] \\
      B \arrow[r] & B_{\lambda} \arrow[r] & B_{\pi}  
    \end{tikzcd}
  \]
  in which the lower horizontal rectangle is cartesian.
  Replacing $g_{\lambda}$ by the projection $E_{\pi} \times_{B_{\pi}} B_{\lambda} \to B_{\lambda}$ and noting that
  \[
    \prolimit_{\lambda \in \Lambda}(E_{\pi} \times_{B_{\pi}} B_{\lambda}) = E_{\pi} \times_{B_{\pi}} (\prolimit_{\lambda \in \Lambda} B_{\lambda}) = E,
  \]
  we may also assume that the lower right square is cartesian.
  Consequently, for every $\lambda \in \Lambda$, also the lower left square is cartesian.
  If $f$ is $\pi$-finite, applying the same reductions to $f$, we may assume that also the upper square is cartesian.
  Hence so is the vertical rectangle, exhibiting the $\pi$-finiteness of $g \circ f$.
  If, one the other hand, $(g \circ f)$ is assumed to be $\pi$-finite, a similar reduction lets us assume that the vertical rectangle is cartesian.
  Hence so is the upper square by cancellation of pullbacks, exhibiting the $\pi$-finiteness of $f\!$.
\end{proof}

Finally, we can prove continuous (un-)straightening.

\begin{theorem}[continuous straightening-unstraightening]
  \label{thm:continuous-un-straightening}
  Let $B$ be a profinite anima and $\Sigma$ a nonempty full class of finite groups.
  \begin{enumerate}[(1)]
  \item\label{thmitem:finite-continuous-un-straightening} 
  The straightening-unstraightening equivalences of the various $\pi$-finite constituents of $B$ assemble into an adjoint equivalence of $\infty$-categories
    \[
      \begin{tikzcd}
        \Fun^{\cts}(B, \Ani_{\Sigma}) \arrow[rr, shift left = 0.2em, "{\Unstr^{B}}"] & & \Ani_{\Sigma}\rel B \arrow[ll, shift left = 0.2em, "{\Str_{B}}"].
      \end{tikzcd}
    \]
    \item \label{thmitem:profinite-continuous-un-straightening}
    If $B$ is in $\Pro(\Ani_{\Sigma})$, there is an induced adjoint equivalence of $\infty$-categories
    \[
      \begin{tikzcd}
        \Pro(\Fun^{\cts}(B, \Ani_{\Sigma})) \arrow[rr, shift left = 0.2em, "{\Unstr^{B}}"] & & \Pro(\Ani_{\Sigma})\rel B \arrow[ll, shift left = 0.2em, "{\Str_{B}}"]
      \end{tikzcd}
    \]
    obtained by passing to pro-categories.
  \end{enumerate}
\end{theorem}

\begin{proof} 
  Choose a level representation $B = \prolimit_{\lambda \in \Lambda} B_{\lambda}$.
  Combining \Cref{lem:straightening-unstraigthening}~\ref{lemitem:pi-finite-(un)-straightening} and~\ref{lemitem:functoriality} with \Cref{lem:filtered-colimits-ff-and-ess-surj}, the straightening-unstraightening equivalences of the various $B_{\lambda}$ assemble into an adjoint equivalence of $\infty$-categories
    \[
      \begin{tikzcd}
        \Fun^{\cts}(B, \Ani_{\Sigma}) \arrow[rr, shift left = 0.2em, "{\Unstr^{B}}"] & & \colim_{\lambda \in \Lambda^{\op}} (\Ani_{\Sigma} \rel B_{\lambda}) \arrow[ll, shift left = 0.2em, "{\Str_{B}}"]
      \end{tikzcd}
    \]
    on the colimits.
    Here, we write $\Unstr^{B} \defeqq \colim_{\lambda \in \Lambda^{\op}}(\Unstr^{B_{\lambda}})$ and $\Str_{B} \defeqq \colim_{\lambda \in \Lambda^{\op}}(\Str_{B_{\lambda}})$ for the induced functors on the colimits.
    \Cref{lem:Sigma-finite-maps}~\ref{lemitem:embedding-to-profinite-anima-over-B} identifies the right-hand side with $\Ani_{\Sigma} \rel B$, showing~\ref{thmitem:finite-continuous-un-straightening}.
    In view of \Cref{lem:Sigma-finite-maps}~\ref{lemitem:equivalence-B-pro-Sigma}, assertion~\ref{thmitem:profinite-continuous-un-straightening} follows from~\ref{thmitem:finite-continuous-un-straightening} by passing to pro-categories.
\end{proof}

\begin{remark}
  \label{rem:continuous-un-straightening-for-pi-finite-anima}
  Let $B$ be a $\pi$-finite anima.
  The functor
  \[
    \Pro(\Fun(B, \Ani_{\pi})) \to \Fun(B, \Pro(\Ani_{\pi}))
  \]
  induced by postcomposing with the inclusion $\Ani_{\pi} \subset \Pro(\Ani_{\pi})$ and then extending along cofiltered limits is an equivalence of $\infty$-categories by \SAG{}{E.6.1.3}.
  Moreover, the diagram
  \[
    \begin{tikzcd}
      \Pro(\Fun(B, \Ani_{\pi})) \arrow[rd, "{\Unstr^{B}}"'] \arrow[r] & \Fun(B, \Pro(\Ani_{\pi})) \arrow[d, "{\colim_{B}}"] \\
      & \Pro(\Ani_{\pi})_{/B}
    \end{tikzcd}
  \]
  commutes.
  Restricted to $\pi$-finite anima, our continuous straightening-unstraightening equivalence therefore coincides with Lurie's variant \cite[\SAGthm{}{E.6.0.7}]{SAG}.
\end{remark}

\begin{corollary}
  \label{cor:mondromy-in-sigma-fully-faithful}
  Let $\Sigma$ be a full class of finite groups and $B$ a profinite anima.
  Precomposition with the pro-$\Sigma$ completion map $u \from B \to B\Spfcompl$ induces a fully faithful embedding
  \[
    u^{\ast} \from \Fun^{\cts}(B\Spfcompl, \Ani_{\Sigma}) \mono \Fun^{\cts}(B, \Ani_{\Sigma}).
  \]
\end{corollary}

\begin{proof}
  Under \Cref{thm:continuous-un-straightening}, the functor $u^{\ast}$ gets identified with
  \[
    \blank \times_{B\Spfcompl} B \from \Ani_{\Sigma}\rel B\Spfcompl \to \Ani_{\Sigma}\rel B.
  \]
  The latter functor is readily seen to be fully faithful by \Cref{rec:Spfcompl-preserves-cofiltered-limits}~\ref{recitem:pfcompl-is-locally-cartesian-and-preserves-limits}.
\end{proof}

\begin{lemma}
  \label{lem:straightening-of-pullback}
  Let
  \[
    \begin{tikzcd}
      E' \arrow[r] \arrow[d, "{f'}"'] \cartesian & E \arrow[d, "{f}"] \\
      B' \arrow[r, "p"'] & B
    \end{tikzcd}
 \]
  be a pullback square in $\Pro(\Ani_{\pi})$ with $f \in \Ani_{\pi} \rel B$.
  Then $\Str_{B'}(f')$ is given by the continuous functor
  \[
    \begin{tikzcd}
      B' \arrow[r, "p"] & B \arrow[r, "{\Str_B(f)}"] & \Ani_{\pi}.
    \end{tikzcd}
 \]
\end{lemma}

\begin{proof}
  Unwinding the definitions, this follows immediately from \Cref{cor:fibrewise-criterion-pullbacks}.
\end{proof}


\section{Applications to profinite shape theory}
\label{sec:profinite-shape-theory}
In this section, we finally prove \Cref{thm:A} of the introduction, see \Cref{thm:lcc-sheaves-and-profinite-shape}, and derive some consequences.

\subsection{The profinite shape}
\label{subsec:profinite-shape}

\begin{definition}
  \label{def:locally-constant-constructible-sheaf}
  Let $\cX$ be an $\infty$-topos with terminal object $1_{\cX}$ and $\Sigma$ a full class of finite groups.
  \begin{enumerate}[(1)]
    \item\label{defitem:lcc-sheaf}
    A sheaf $F \in \cX$ is called \emph{locally $\Sigma$-constant constructible} if there exists an effective epimorphism $\sqcup_{i = 1}^n U_i \surj 1_{\cX}$ and $\Sigma$-finite anima $K_1, \ldots, K_n$ together with equivalences $F \times U_i \simeq \Gamma^{\ast}K_i \times U_i$ in $\cX_{\!/U_i}$.
    \item \label{defitem:subcategory-of-Sigma-lcc-sheaves}
    We write $\lcc(\cX, \Ani_{\Sigma}) \subset \cX$ for the full subcategory spanned by the locally $\Sigma$-constant constructible sheaves.

    \item \label{defitem:subcategory-of-lcc-sheaves}
    When $\Sigma = \pi$, we simply say \emph{locally constant constructible} and write $\lcc(\cX) \defeqq \lcc(\cX, \Ani_{\pi})$.

    \item \label{defitem:lcc-maps}
    A morphism $f \from F' \to F$ in $\cX$ is called \emph{locally constant constructible} if $f \in \lcc(\cX_{\!/F})$.
  \end{enumerate}
\end{definition}

\begin{recollection}
\label{rec:closure-properties-of-lcc}
Let $f_{\!\ast} \from \cX \to \cS$ be a geometric morphism of $\infty$-topoi.
\begin{enumerate}[(1)]
    	\item \label{recitem:closed} 
    	The subcategory $\lcc(\cX)\subseteq \cX$ is closed under finite limits and finite coproducts, see \SAG{}{E.2.5.2} and \SAG{}{E.2.6.8}.

   	 \item \label{recitem:pullback-preserves-lcc} 
    	The pullback functor $f^{\ast} \from \cS \to \cX$ preserves being locally $\Sigma$-constant constructible. 
\end{enumerate}
\end{recollection}

\begin{lemma}
  \label{lem:profinite-shape}
  Let $f_{\!\ast} \from \cX \to \cS$ and $g_{\ast} \from \cY \to \cX$ be geometric morphisms of $\infty$-topoi.
  \begin{enumerate}[(1)]
    \item \label{lemitem:adjoints} 
    Each functor in the commutative diagram
    \[
      \begin{tikzcd}
        \cS \arrow[r, "{f^{\ast}}"] & \cX \\
        \lcc(\cS) \arrow[u, hook] \arrow[r, "{f^{\ast}}"'] & \lcc(\cX) \arrow[u, hook]
      \end{tikzcd}
    \]
    has a pro-left adjoint and the resulting diagram
    \[
      \begin{tikzcd}
        \Pro(\cX) \arrow[r, "{f_{\! !}}"] \arrow[d, "{(\blank)\lcccompl}"'] & \Pro(\cS) \arrow[d, "{(\blank)\lcccompl}"] \\
        \Pro(\lcc(\cX)) \arrow[r, "{f_{\!\sharp}}"'] & \Pro(\lcc(\cS))
      \end{tikzcd}
    \]
    commutes.
    \item \label{lemitem:shape-composition} There is a natural equivalence $(f \circ g)_{\sharp} \simeq f_{\!\sharp} \circ g_{\sharp}$.
  \end{enumerate}
\end{lemma}

\begin{proof}
  Assertion~\ref{lemitem:adjoints} is a formal consequence of \Cref{rec:closure-properties-of-lcc}~\ref{recitem:closed} and~\ref{recitem:pullback-preserves-lcc}.
  Assertion~\ref{lemitem:shape-composition} follows from $(f \circ g)^{\ast} \simeq g^{\ast} \circ f^{\ast}$ by passing to left adjoints.
\end{proof}

\begin{definition}[profinite shape]
  \label{def:profinite-shape}
  Let $f = (f^{\ast} \ladj f_{\!\ast}) \from \cX \to \cS$ be a geometric morphism of $\infty$-topoi.
  \begin{enumerate}[(1)]
  \item\label{defitem:relative} 
    The functor
    \[
       f_{\!\sharp} \from \Pro(\lcc(\cX)) \to \Pro(\lcc(\cS)) 
    \] 
    is called the \emph{(relative) profinite shape functor of $f$}.
    The sheaf $f_{\!\sharp}F$ is called the \emph{profinite shape of $F$ relative to $f$}.
    \item \label{defitem:shape-of-morphism}
    The sheaf
    \[
      \pfhtpytype(\cX/\cS) \defeqq \pfhtpytype(f) \defeqq f_{\!\sharp} 1_{\cX}
    \]
    is called the \emph{relative profinite shape of $f$}.
    \item \label{defitem:absolute}
    The functor
    \[
      \pfhtpytype \defeqq \Gamma_{\cX, \sharp} \from \Pro(\lcc(\cX)) \to \Pro(\Ani_{\pi})
    \]
    is called the \emph{profinite shape functor of $\cX$}.
    The profinite anima $\pfhtpytype(F)$ is called the \emph{profinite shape} of $F \in \Pro(\lcc(\cX))$.

    \item \label{defitem:shape-of-topos}
    The profinite anima
    \[
       \pfhtpytype(\cX) \defeqq \pfhtpytype(1_{\cX})
     \]
     is called the \emph{profinite shape of $\cX$}.
  \end{enumerate}
\end{definition}

\begin{remark}
  Let $f_{\!\ast} \from \cX \to \cS$ be geometric morphism of $\infty$-topoi and $\Gamma_{\ast} \from \cS \to \Ani$ the unique geometric morphism.
  \begin{enumerate}[(1)]
    \item An anima $K$ is a locally constant constructible object of $\Ani$ if and only if it is $\pi$-finite.
      Consequently, lcc-completion $(\blank)\lcccompl$ on $\Pro(\Ani)$ agrees with profinite completion $(\blank)\pfcompl$.
    \item By \Cref{lem:profinite-shape}, we see that
      \[
        \pfhtpytype(\cX) = (\Gamma_{!}1_{\cX})\pfcompl
      \]
      recovers the usual definition of the profinite shape.
      \item By definition, $\pfhtpytype(\cX/\cS)$ has the following universal property:
        Functorially in $F \in \lcc(\cS)$, there is an equivalence
      \[
         \Map_{\Pro(\lcc(\cS))}(\pfhtpytype(\cX/\cS), F) \simeq \Gamma_{\ast}(\cX, f^{\ast}F).
      \]
      That is: $\pfhtpytype(\cX/\cS)$ encodes the cohomology of locally constant constructible sheaves that arise via pullback from $\cS$.
  \end{enumerate}
\end{remark}

\begin{lemma}
  \label{lem:restricted-shape-adjunction}
  Let $\cX$ be an $\infty$-topos.
  \begin{enumerate}[(1)]
    \item \label{lemitem:induced-adjunction} The adjunction
    $\begin{tikzcd}
      \pfhtpytype \from \Pro(\lcc(\cX)) \arrow[r, shift left = 0.2em] & \Pro(\Ani_{\pi}) \from \Gamma^{\ast} \arrow[l, shift left = 0.2em]
    \end{tikzcd}$ induces an adjunction
    \[
      \begin{tikzcd}
        \relpfhtpytype{\cX} \from \Pro(\lcc(\cX)) \arrow[rr, shift left = 0.2em] && \Pro(\Ani_{\pi})_{/\pfhtpytype(\cX)} \from \Gamma_{\pfhtpytype(\cX)}^{\ast}. \arrow[ll, shift left = 0.2em]
      \end{tikzcd}
    \]
    Here, $\Gamma_{\pfhtpytype(\cX)}^{\ast}$ is given as the composition of the sliced pullback functor $\Gamma^{\ast} \from \Pro(\Ani_{\pi})_{/\pfhtpytype(\cX)} \to \Pro(\lcc(\cX))_{/\Gamma^{\ast}\!\pfhtpytype(\cX)}$ followed by pulling back along the canonical section $1_{\cX} \to \Gamma^{\ast}\!\pfhtpytype(\cX)$ given by the unit of the adjunction $\pfhtpytype \ladj \Gamma^{\ast}$.
    \item \label{lemitem:restricted-right-adjoint}
    The functor $\Gamma^{\ast}_{\pfhtpytype(\cX)}$ restricts to a functor $\Gamma^{\ast}_{\pfhtpytype(\cX)} \from \Ani_{\pi} \rel \pfhtpytype(\cX) \to \lcc(\cX)$.
  \end{enumerate}
\end{lemma}

\begin{proof}
  Assertion \lref{lemitem:induced-adjunction} is formal, and \lref{lemitem:restricted-right-adjoint} is a direct consequence of \Cref{rec:closure-properties-of-lcc}~\ref{recitem:closed}.
\end{proof}

\subsection{The monodromy equivalence}
\label{subsec:monodromy-equivalence}

For a geometrically pointed connected qcqs scheme $(X, x)$, the category of locally constant constructible \'etale sheaves of sets  on $X$ is equivalent to finite sets with a continuous action by the \'etale fundamental group $\etfdtlgrp(X, x)$.
Under this equivalence, a sheaf $F$ is mapped to the monodromy action of $\etfdtlgrp(X, x)$ on the stalk $F_{\!x}$.
More canonically, and without connectivity assumptions on $X$, the locally constant constructible sheaf $F$ can be recovered from a continuous local system on $\truncethtpytype{\leq 1}(X)$, the \emph{monodromy representation} $\Mon(F) \from \truncethtpytype{\leq 1}(X) \to \Set_{\fin}$ of the \'etale fundamental groupoid of $X$.
There is an equivalence of $1$-categories
\[
    \Mon \from \lcc(X_{\et}, \Set) \longrightarrow \Fun^{\cts}(\truncethtpytype{\leq 1}(X), \Set_{\fin}).
\]
The fundamental theorem of profinite shape theory, due to Lurie \cite[\SAGapp{E}]{SAG}, Bachmann--Hoyois \cite[Prop. 10.1]{bachmann2020} and Barwick--Glasman--Haine \cite[Prop. 4.4.18]{barwick2020}, refines the above equivalence to one capturing the whole of $\lcc(X_{\et})$.

\begin{theorem}[monodromy equivalence]
\label{thm:monodromy-equivalence}
Let $\cX$ be an $\infty$-topos.
The composition $\Mon = \Str_{\pfhtpytype(\cX)} \circ \relpfhtpytype{\cX}$ is an equivalence of $\infty$-categories
\[
	\Mon  \from \lcc(\cX) \longrightarrow \Fun^{\cts}(\pfhtpytype(\cX), \Ani_{\pi})
\]
that carries a locally constant constructible sheaf $F$ to its \emph{monodromy representation} $\Mon(F)$.
Its inverse is given by the composition
\[
	\begin{tikzcd}
		\Fun^{\cts}(\pfhtpytype(\cX), \Ani_{\pi}) \arrow[r, "{\Unstr^{\pfhtpytype(\cX)}}"] 
		& \Ani_{\pi} \rel \pfhtpytype(\cX) \arrow[r, "{\Gamma^{\ast}_{\pfhtpytype(\cX)}}"] 
		& \lcc(\cX).
	\end{tikzcd}
\]
\end{theorem}

\begin{proof}
We refer the reader to \cite[Prop. 4.4.18]{barwick2020} and \cite[Prop. 10.1]{bachmann2020} (together with \cite[Thm. 4.3]{hoyois2017}).
Note that both references just state the existence of ``canonical'' equivalences.
However, upon inspecting the proofs in either of the above cases, one sees the explicit description we stated.
\end{proof}

\begin{remark}
  \label{rem:relation-with-classical-monodromy-equivalence}
  Let $\cX$ be an $\infty$-topos.
  \begin{enumerate}[(1)]
    \item\label{remitem:classical-monodromy}
    In view of \Cref{thm:relation-to-classical-local-systems}, the monodromy equivalence restricts to an equivalence
    \[
      \Fun^{\cts}(\pfhtpytype(\cX), \Set_{\fin}) \isomto \lcc(\cX, \Set)
    \]
     that in turn induces equivalences
    \[
      \Fun^{\cts}(\pfhtpytype(\cX), \Grp_{\fin}) \isomto \lcc(\cX, \Grp_{\fin}) \quad \text{and} \quad \Fun^{\cts}(\pfhtpytype(\cX), \Ab_{\fin}) \isomto \lcc(\cX, \Ab_{\fin}).
    \]
    \item \label{remitem:classical-etale-monodromy}
    If $X$ is a connected scheme, the choice of a geometric point $\bar x \to X$ determines an equivalence
    \[
      \Fun^{\cts}(\ethtpytype(X), \Set_{\fin}) = \Fun^{\cts}(\truncethtpytype{\leq 1}(X), \Set_{\fin}) \simeq \etfdtlgrp(X, \bar x)\textnormal{-}\!\Set_{\fin},
    \]
    where the right-hand side denotes the category of finite sets equipped with a continuous $\etfdtlgrp(X, \bar x)$-action.
    Under this equivalence, \ref{remitem:classical-monodromy} recovers \stacks{0DV5}.
  \end{enumerate}

\end{remark}

In view of \Cref{cor:mondromy-in-sigma-fully-faithful}, we are led to the following notion.

\begin{definition}[monodromy in $\Sigma$]
  \label{def:sheaf-with-monodromy-in-Sigma}
  Let $\cX$ be an $\infty$-topos, $\Sigma$ a full class of finite groups, and $F \in \lcc(\cX, \Ani_{\Sigma})$.
  \begin{enumerate}[(1)]
    \item\label{defitem:monodromy-in-Sigma}
    The sheaf $F$ has \emph{monodromy in $\Sigma$} if its monodromy representation $\Mon(F)$ lies in the essential image of the inclusion $\Fun^{\cts}(\pfhtpytype(\cX)\Spfcompl, \Ani_{\Sigma}) \mono \Fun^{\cts}(\pfhtpytype(\cX), \Ani_{\Sigma})$ of \Cref{cor:mondromy-in-sigma-fully-faithful}.
    \item \label{defitem:subcategory-of-monodromy-in-Sigma}
    We write $\lcc_{\Sigma}(\cX)$ for the full subcategory of $\lcc(\cX, \Ani_{\Sigma})$ spanned by the locally $\Sigma$-constant constructible sheaves with monodromy in $\Sigma$.
  \end{enumerate}
\end{definition}

\subsection{The profinite shape equivalence}
\label{subsec:shape-equivalence}

The monodromy equivalence has the following fairly straightforward consequence, that appears to not have been recorded in the literature so far.

\begin{theorem}[profinite shape equivalence]
\label{thm:lcc-sheaves-and-profinite-shape}
Let $\cX$ be an $\infty$-topos.
\begin{enumerate}[(1)]
	\item \label{thmitem:shape-equivalence} 
        The profinite shape adjunction 
	\[
      		\begin{tikzcd}
        			\Pro(\lcc(\cX)) \arrow[rr, shift left = 0.2em, "{\relpfhtpytype{\cX}}"] 
			&& \Pro(\Ani_{\pi})_{\!/\!\pfhtpytype(\cX)}
			\arrow[ll, shift left = 0.2em, "{\Gamma_{\pfhtpytype(\cX)}^{\ast}}"]
      		\end{tikzcd}
    	\]
        is an adjoint equivalence of $\infty$-categories.

    	\item \label{thmitem:commutes}
    	Let $f = (f^{\ast}, f_{\! \ast}) \from \cX \to \cS$ be a geometric morphism of $\infty$-topoi. 
	In the diagram
    	\[
      		\begin{tikzcd}
        		\Pro(\lcc(\cS)) \arrow[rr, shift right = 0.2em, "{f^{\ast}}"'] \arrow[dd, shift right = 0.2em, "{\relpfhtpytype{\cS}}"'] 
		&& \Pro(\lcc(\cX)) \arrow[ll, shift right = 0.2em, "{f_{\sharp}}"'] 
		\arrow[dd, shift right = 0.2em, "{\relpfhtpytype{\cX}}"'] 
		\\
        		&& \\
        		\Pro(\Ani_{\pi})_{\!/\!\pfhtpytype(\cS)} 
		\arrow[uu, shift right = 0.2em, "{\Gamma^{\ast}_{\pfhtpytype(\cS)}}"'] 
		\arrow[rr, shift right = 0.2em, "{\blank \times_{\pfhtpytype(\cS)} \pfhtpytype(\cX)}"'] 
		&& \Pro(\Ani_{\pi})_{\!/\!\pfhtpytype(\cX)} 
		\arrow[ll, shift right = 0.2em, "{\pfhtpytype(f) \circ \blank}"'] 
		\arrow[uu, shift right = 0.2em, "{\Gamma^{\ast}_{\pfhtpytype(\cX)}}"']
		\\
      		\end{tikzcd}
    	\]
    	both the subdiagrams of left-adjoints and of right-adjoints commute.
\end{enumerate}
\end{theorem}

\begin{proof}
By \Cref{thm:continuous-un-straightening}~\ref{thmitem:finite-continuous-un-straightening}, the unstraightening functor $\Unstr^{\pfhtpytype(\cS)}$ is an equivalence.
Therefore \Cref{thm:monodromy-equivalence} implies that also $\Gamma^{\ast}_{\pfhtpytype(\cS)} \from \Ani_{\pi} \rel \pfhtpytype(\cS) \to \lcc(\cS)$ is an equivalence.
Passing to pro-categories, we see that also
\[
	\Gamma^{\ast}_{\pfhtpytype(\cS)} \from \Pro(\Ani_{\pi})_{\!/\!\pfhtpytype(\cS)} \to \Pro(\lcc(\cS))
\]
is an equivalence by \Cref{lem:Sigma-finite-maps}~\ref{lemitem:equivalence-B-pro-Sigma}.
By construction, $\relpfhtpytype{\cS}$ is left-adjoint to $\Gamma^{\ast}_{\pfhtpytype(\cS)}$, hence also an equivalence.
This shows~\ref{thmitem:shape-equivalence}.
That the subdiagram of right-adjoints in \ref{thmitem:commutes} commutes is obvious.
Hence so does the one of left-adjoints by uniqueness of adjoints.
\end{proof}

\begin{remark}
  \label{rem:relation-to-SAG}
  Alternatively, it should be possible to deduce \Cref{thm:lcc-sheaves-and-profinite-shape} by carefully tracing the proof of \SAG{}{A.8.0.5} in the case $\cC = \lcc(\cX)$.
\end{remark}

\begin{corollary}
  \label{cor:restricted-shape-equivalence}
  Let $\cX$ be an $\infty$-topos.
  The profinite shape equivalence restricts to an adjoint equivalence of $\infty$-categories
      $\begin{tikzcd}
        \relpfhtpytype{\cX} \from \lcc(\cX) \arrow[rr, shift left = 0.2em] && \Ani_{\pi}\rel\pfhtpytype(\cX)\from \Gamma_{\pfhtpytype(\cX)}^{\ast}. \arrow[ll, shift left = 0.2em]
      \end{tikzcd}$
\end{corollary}

\begin{proof}
Since both $\lcc(\cX)$ and $\Ani_{\pi} \rel \pfhtpytype(\cX)$ are essentially small and idempotent complete, this follows from \Cref{thm:lcc-sheaves-and-profinite-shape}~\ref{thmitem:shape-equivalence} and \Cref{lem:Sigma-finite-maps}~\ref{lemitem:equivalence-B-pro-Sigma} by passing to cocompact objects, see \HTT{}{5.4.2.4} and \HTT{}{5.4.3.6}.
\end{proof}

\begin{corollary}
  \label{cor:monodromy-characterisation-of-lcc-sheaves}
  Let $\cX$ be an $\infty$-topos and $f \from F \to G$ as well as $g \from G \to H$ two morphisms in $\cX$.
  \begin{enumerate}[(1)]
  \item\label{coritem:lcc-is-pullback} 
  If $F \in \lcc(\cX)$, the canonical diagram
  \[
    \begin{tikzcd}
      F \arrow[r] \arrow[d] \cartesian & \Gamma^{\ast}\!\pfhtpytype(F) \arrow[d] \\
      1_{\cX} \arrow[r] & \Gamma^{\ast}\!\pfhtpytype(\cX)
    \end{tikzcd}
  \]
  is a pullback in $\Pro(\lcc(\cX))$.

  \item \label{coritem:lcc-maps-stable-under-composition}
  If both $f$ and $g$ are locally constant constructible, so is $g \circ f$.
  \end{enumerate}
\end{corollary}

\begin{proof}
  If $F \in \lcc(\cX)$, we have $F \simeq \Gamma^{\ast}_{\pfhtpytype(\cX)}\!\pfhtpytype(F)$.
  Unwinding the definition of $\Gamma^{\ast}_{\pfhtpytype(\cX)}$, this precisely gives the pullback square of~\ref{coritem:lcc-is-pullback}.
  To see~\ref{coritem:lcc-maps-stable-under-composition}, note that we may assume $H = 1_{\cX}$ by replacing $\cX$ with $\cX_{\!/H}$.
  Then by~\ref{coritem:lcc-is-pullback}, the squares
  \[
    \begin{tikzcd}
      F \arrow[d, "{(u_F, f)}"'] \arrow[rr, "{f}"] && G \arrow[d, "{(u_G, \id_{G})}"] & G \arrow[rr, "{g}"] \arrow[d, "{u_G}"'] && 1_{\cX} \arrow[d, "{u_{\cX}}"] \\
      \Gamma^{\ast}\!\pfhtpytype(F) \times G \arrow[rr, "{\Gamma^{\ast}\!\pfhtpytype(f) \times \id_G}"'] && \Gamma^{\ast}\!\pfhtpytype(G) \times G & \Gamma^{\ast}\!\pfhtpytype(G) \arrow[rr, "{\Gamma^{\ast}\!\pfhtpytype(g)}"'] && \Gamma^{\ast}\!\pfhtpytype(\cX).
    \end{tikzcd}
  \]
  are cartesian in $\Pro(\lcc(\cX_{\!/G}))$ respectively $\Pro(\lcc(\cX))$, hence both are cartesian in $\Pro(\cX)$.
  It follows that also
  \[
    \begin{tikzcd}
      F \arrow[r, "{u_{F}}"] \arrow[d, "{g \circ f}"'] \cartesian & \Gamma^{\ast}\!\pfhtpytype(F) \arrow[d, "{\Gamma^{\ast}\!\pfhtpytype(g \circ f)}"] \\
      1_{\cX} \arrow[r, "{u_{H}}"'] & \Gamma^{\ast}\!\pfhtpytype(\cX)
    \end{tikzcd}
  \]
  is cartesian in $\Pro(\cX)$.
  We conclude by combining \Cref{cor:pi-finite-maps-stable-under-composition}, \Cref{rec:closure-properties-of-lcc}~\ref{recitem:closed}, and \Cref{cor:restricted-shape-equivalence}.
\end{proof}

\begin{corollary}
  \label{cor:shape-of-pullback}
  Let $f = (f^{\ast}, f_{\!\ast}) \from \cX \to \cS$ be a geometric morphism of $\infty$-topoi, $\Sigma$ a full class of finite groups, and $F \in \lcc(\cS)$.
  \begin{enumerate}[(1)]
  \item\label{coritem:pullback-square} 
  The canonical diagram
  \[
    \begin{tikzcd}
      \pfhtpytype(f^{\ast}F) \arrow[r] \arrow[d] \cartesian & \pfhtpytype(F) \arrow[d] \\
      \pfhtpytype(\cX) \arrow[r, "{\pfhtpytype(f)}"'] & \pfhtpytype(\cS)
    \end{tikzcd}
  \]
  is a pullback square in $\Pro(\Ani_{\pi})$.

  \item \label{coritem:monodromy-of-pullback}
  The monodromy representation of $f^{\ast} F$ is given as the continuous functor
  \[
    \begin{tikzcd}
      \pfhtpytype(\cX) \arrow[r, "{\pfhtpytype(f)}"] & \pfhtpytype(\cS) \arrow[r, "{\Mon(F)}"] & \Ani_{\pi}.
    \end{tikzcd}
  \]

  \item \label{coritem:monodromy-in-Sigma-stable-under-pullback}
  If $F \in \lcc_{\Sigma}(\cS)$, then $f^{\ast}F \in \lcc_{\Sigma}(\cX)$.
  \end{enumerate}
\end{corollary}

\begin{proof}
  Assertion~\ref{coritem:pullback-square} is a formal consequence of \Cref{thm:lcc-sheaves-and-profinite-shape}, since the vertical adjunctions are equivalences:
  Since $F = \Gamma^{\ast}_{\pfhtpytype(\cS)}\!\pfhtpytype(F)$, we have that $f^{\ast} F = \Gamma^{\ast}_{\pfhtpytype(\cX)} (\pfhtpytype(F) \times_{\pfhtpytype(\cS)} \pfhtpytype(\cX))$ by \Cref{thm:lcc-sheaves-and-profinite-shape}~\ref{thmitem:commutes}.
  The result follows by applying $\pfhtpytype(\blank)$.
  Assertion~\ref{coritem:monodromy-of-pullback} follows from~\ref{coritem:pullback-square} and \Cref{lem:straightening-of-pullback}.
  Finally,~\ref{coritem:monodromy-in-Sigma-stable-under-pullback} follows from~\ref{coritem:monodromy-of-pullback} by considering the commutative diagram of continuous functors
  \[
    \begin{tikzcd}
      \pfhtpytype(\cX) \arrow[rr, "{\pfhtpytype(f)}"] \arrow[d] && \pfhtpytype(\cS) \arrow[d] \arrow[rr, "{\Mon(F)}"] && \Ani_{\Sigma}. \\
      \pfhtpytype(\cX)\Spfcompl \arrow[rr, "{\pfhtpytype(f)\Spfcompl}"'] && \pfhtpytype(\cS)\Spfcompl \arrow[rru, "{{\exists?}}"', dashed] & & 
    \end{tikzcd}
    \qedhere
  \]
\end{proof}

Finally, we also have the following pro-$\Sigma$ variant of the profinite shape equivalence:

\begin{theorem}[pro-$\Sigma$ shape equivalence]
  \label{thm:pro-Sigma-shape-equivalence}
  Let $\cX$ be an $\infty$-topos and $\Sigma$ full class of finite groups.
  The pro-$\Sigma$ shape induces an adjoint equivalence of $\infty$-categories
  \[
      \begin{tikzcd}
        \pfhtpytype(\blank)\Spfcompl \from \Pro(\lcc_{\Sigma}(\cX)) \arrow[rr, shift left = 0.2em]
        && \Pro(\Ani_{\Sigma})_{\!/\!\pfhtpytype(\cX)\Spfcompl} \from \Gamma_{\pfhtpytype(\cX)\Spfcompl}^{\ast},
        \arrow[ll, shift left = 0.2em]
      \end{tikzcd}
  \]
  where $\Gamma_{\pfhtpytype(\cX)\Spfcompl}^{\ast} = \Gamma^{\ast}_{\pfhtpytype(\cX)} \circ\ (\blank \times_{\pfhtpytype(\cX)\Spfcompl} \pfhtpytype(\cX))$.
\end{theorem}

\begin{proof}
  Combining \Cref{rec:closure-properties-of-lcc}~\ref{recitem:closed} with \Cref{cor:shape-of-pullback}~\ref{coritem:monodromy-in-Sigma-stable-under-pullback} and observing that $\Ani_{\Sigma} = \lcc_{\Sigma}(\Ani) \subset \lcc(\Ani) = \Ani_{\pi}$, we see that all of the functors in the commutative diagram on the left
  \[
    \begin{tikzcd}
      \lcc(\cX) \arrow[rr, leftarrow, "{\Gamma^{\ast}}"] && \Ani_{\pi} & \Pro(\lcc(\cX)) \arrow[dd, shift right = 0.3em, "{(\blank)\Spfcompl}"'] \arrow[rr, shift left = 0.3em, "{\pfhtpytype}"] && \Pro(\Ani_{\pi})_{/\pfhtpytype(\cX)} \arrow[ll, shift left = 0.3em, "{\Gamma^{\ast}_{\pfhtpytype(\cX)}}"] \arrow[dd, shift right = 0.3em, "{(\blank)\Spfcompl}"']\\
      && & && \\
      \lcc_{\Sigma}(\cX) \arrow[uu, hook] \arrow[rr, leftarrow, "{\Gamma^{\ast}}"'] && \Ani_{\Sigma} \arrow[uu, hook] & \Pro(\lcc_{\Sigma}(\cX)) \arrow[uu, hook, shift right = 0.3em] \arrow[rr, shift left = 0.3em, "{\Gamma_{\sharp, \Sigma}}"] && \Pro(\Ani_{\Sigma})_{/\pfhtpytype(\cX)\Spfcompl} \arrow[ll, shift left = 0.3em, "{\Gamma^{\ast}_{\pfhtpytype(\cX)\Spfcompl}}"] \arrow[uu, shift right = 0.3em, hook, "{\blank \times_{\pfhtpytype(\cX)\Spfcompl} \pfhtpytype(\cX)}"']
    \end{tikzcd}
  \]
  admit pro-left adjoints, inducing the commutative diagram of sliced adjoints on the right.
  By \Cref{thm:lcc-sheaves-and-profinite-shape}, the top horizontal adjunction is an equivalence.
  Unwinding the definition of $\lcc_{\Sigma}(\cX)$, we moreover see that the image of $\Pro(\lcc_{\Sigma}(\cX))$ under $\pfhtpytype$ coincides with the image of the pullback functor $(\blank \times_{\pfhtpytype(\cX)\Spfcompl} \pfhtpytype(\cX))$.
  Since $\Gamma_{\sharp, \Sigma} = (\blank)\Spfcompl\ \circ \pfhtpytype$, the claim follows.
\end{proof}


\printbibliography

\end{document}
